\def\E{\end{document}}
\documentclass[10pt]{article}
\usepackage{mathrsfs}
\usepackage{amssymb}
\usepackage{amssymb,amsmath,color}

\usepackage{amsmath,amsthm,bm,mathtools}
\newcommand{\T}{{\mathbb T_L}}

\newcommand{\Div}{\operatorname{div}}
\newcommand{\supp}{\operatorname{supp}}
\newcommand{\rb}{\bar\rho}
\newcommand{\ub}{\bar{\bm u}}
\newcommand{\bb}{\bar{\bm b}}
\newcommand{\tb}{\bar\theta}
\newcommand{\Pb}{\bar P}
\newcommand{\Prb}{\bar P_\rho}
\newcommand{\Ptb}{\bar P_\theta}
\newcommand{\norm}[1]{\lVert #1\rVert}
\newcommand{\dual}[2]{\langle #1,#2\rangle}

\begin{document}
\title{\bf Local exact boundary controllability for the full compressible magnetohydrodynamic equations}

\author{Qiang Tao
\thanks {School of Mathematical Sciences, Shenzhen University, 518060, Shenzhen, China
(taoq060@126.com, taoq@szu.edu.cn).}
}

\date{}
\maketitle
\newtheorem{theorem}{Theorem}[section]
\newtheorem{definition}{Definition}[section]
\newtheorem{lemma}[theorem]{Lemma}
\newtheorem{geometriclemma}{Lemma}
\renewcommand{\thegeometriclemma}{3.1a}
\newtheorem{proposition}[theorem]{Proposition}
\newtheorem{corollary}{Corollary}[section]
\newtheorem{remark}{Remark}[section]
\renewcommand{\theequation}{\thesection.\arabic{equation}}
\catcode`@=11 \@addtoreset{equation}{section} \catcode`@=12

{\bf Abstract.}\ \ {\small In this paper, we are concerned with the local exact boundary controllability for the full compressible magnetohydrodynamic equations in a bounded domain in two and three space dimensions. We establish the local exact controllability around a constant equilibrium state. Our proof relies on the observability and controllability analysis of a suitably reformulated linearized system, together with a fixed point argument in weighted Sobolev spaces. The main novelty of this paper is to overcome the structural difficulties arising from the divergence free constraint and the strong coupling between the fluid and magnetic variables. To accommodate the solenoidal constraint, we introduce a magnetic potential formulation and develop a mixed coordinate approach.}

{\bf Mathematics Subject Classification.}\ \  {93B05, 76W05, 76X05}

{\bf Keywords.}\ \ {\small compressible magnetohydrodynamic equations; nonisentropic flows; controllability}

\section{Introduction}
Magnetohydrodynamics (MHD) studies the dynamics of electrically conducting fluids
and the theory of the macroscopic interaction of electrically conducting fluids with a
magnetic field. The governing equations for compressible nonisentropic MHD
flows in $\mathbb{R}^d$  with $d = 2,3$, derived from fluid mechanics with appropriate modifications to account for
electrical forces, have the following form in the Eulerian coordinates:
\begin{equation}\label{1.1}
\left\{\begin{aligned}
 &\partial_t\widetilde\rho+\Div(\widetilde\rho\bm U)=0,\\
 &\widetilde\rho(\partial_t\bm U+\bm U\cdot\nabla\bm U)
 -\mu\Delta\bm U-(\lambda'+\mu)\nabla\Div\bm U
 +\nabla P(\widetilde\rho,\widetilde\theta)=(\nabla\times\bm B)\times\bm B,\\
 &c_v\widetilde\rho(\partial_t\widetilde\theta+\bm U\cdot\nabla\widetilde\theta)
 -\kappa\Delta\widetilde\theta+P(\widetilde\rho,\widetilde\theta)\Div\bm U
 =2\mu|D(\bm U)|^2+\lambda'(\Div\bm U)^2+\nu|\nabla\times\bm B|^2,\\
 &\partial_t\bm B-\nabla\times(\bm U\times\bm B)-\nu\Delta\bm B=0,
 \qquad\Div\bm B=0,
\end{aligned}\right.
\end{equation}
where $\widetilde\rho$ denotes the density, $\bm U $ the velocity, $\widetilde{\theta}$ the temperature, $\bm B$
the magnetic field and pressure $P=P(\widetilde{\rho}, \widetilde{\theta})$ is a smooth function of the density and the temperature.  
The viscosity coefficients $\mu$ and $\lambda'$ satisfy the physical restrictions $\mu>0$ and $2\mu+d\lambda'\geq0$. 
In particular, $\lambda:=2\mu+\lambda'\ge2\mu(1-1/d)>0$ for $d=2,3$.
The positive constants $c_v$ and $\kappa$ are the specific heat at constant volume and the heat conductivity coefficient, and the constant $\nu>0$ is the resistivity coefficient which is inversely proportional to the
electrical conductivity constant and acts as the magnetic diffusivity of magnetic fields. The symbol $\otimes$ denotes the
tensor product. $D(\bm U)$ is the
deformation tensor given by
$$D(\bm U )=\frac{1}{2}[\nabla \bm U +(\nabla\bm U )^{T}].$$
We use $\nabla\times$ for the usual curl in three dimensions. In two dimensions,
\[
 \nabla\times a=(\partial_2a,-\partial_1a),\qquad
 \nabla\times\bm v=\partial_1v_2-\partial_2v_1,
 \qquad \bm U\times\bm B=U_1B_2-U_2B_1.
\]
For a scalar $a$, the planar products are $a\times\bm v=(-av_2,av_1)$ and $\bm v\times a=(v_2a,-v_1a)$, as obtained by embedding planar vectors in $\mathbb R^3$.

The application of MHD
finds a very wide range of physical objects, from liquid metals to
cosmic plasmas. We refer to \cite{C1970, KL1965} for more introduction to the compressible MHD equations.
However, the problems of well-posedness and dynamical behaviors of MHD systems are rather complicated to investigate because both hydrodynamic and electrodynamic effects have to be considered.
Although the electric field $E$ does not appear in the systems, it is indeed induced by moving
the conductive flow in the magnetic field, according to the following relation:
$$E=\nu\nabla\times\mathbf{B}-\bm U \times\mathbf{B}.$$

Due to its physical relevance and the analytical difficulties arising from the coupling of transport and dissipative mechanisms, the controllability of compressible flows has attracted considerable attention. For the compressible Euler equations, in \cite{LR2003} the local exact controllability for the
33 one-dimensional compressible Euler equations was studied in the context of classical solutions. Glass \cite{G2007,G2014} subsequently established exact controllability results for the one-dimensional isentropic and nonisentropic Euler equations in the framework of entropy solutions. In \cite{N2011} Nersisyan proved approximate controllability of the three-dimensional compressible Euler system by means of a finite-dimensional control acting only on the velocity equation. For viscous compressible flows, 
In \cite{A2011}, Amosova considered the controllability for one-dimensional compressible viscous fluid when the equations are stated in Lagrangian forms and the density coincides with the target density. Based on a fixed point argument on a kind of linearized system,
Ervedoza et al. \cite{EGGP2012} established the local exact controllability to constant states for the one-dimensional isentropic Navier-Stokes equations. Then, in \cite{EGG2016}, they proved the multidimensional isentropic flows is local exact boundary controllable by introducing a new variable similar to the so-called viscous effective flux in the adjoint system as well. Later, these results had been extended to the nonisentropic case in \cite{M2019}. Ervedoza and Savel \cite{ES2018} further obtained local boundary controllability around smooth nonconstant trajectories in one space dimension under a geometric condition on the target flow. More recently, Koike et al. \cite{KSV2025} established a local exact Lagrangian controllability result for the one-dimensional barotropic compressible Navier--Stokes equations. For controllability properties of linearized compressible Navier--Stokes systems, we refer in particular to \cite{CM2015,CMRR2014,CRR2012,Maity2015}. 

Controllability results for magnetohydrodynamic flows are considerably more limited, especially in the compressible setting. For incompressible viscous MHD systems, pioneering local exact controllability results were obtained by Barbu and collaborators \cite{BHPS2003,BHPS2005,HPS2006,HPS2007}, and Badra \cite{B2014} subsequently established local controllability to trajectories by localized internal controls acting on both the velocity and the magnetic field. In \cite{G2013} Galan proved approximate controllability of the incompressible MHD equations on the three-dimensional torus by the return method. More recently, substantial progress has been made for incompressible MHD. Rissel and Wang \cite{RW2021} showed global exact boundary controllability of ideal incompressible MHD flows in a planar duct, while Kukavica, Novack, and Vicol \cite{KNV2022} obtained an exact boundary controllability result for ideal MHD in rectangular domains and identified the associated compatibility conditions. For viscous incompressible MHD, Rissel and Wang \cite{RW2024} established small-time global approximate controllability under coupled Navier slip boundary conditions. Further developments include exact controllability of two-dimensional ideal incompressible MHD \cite{R2025} and local exact controllability to trajectories for the two-dimensional viscous incompressible MHD system using distributed controls acting only on the velocity equation \cite{TYY2025}. In contrast, much less is known for compressible MHD. In \cite{T2018}, local exact boundary controllability was established for the planar compressible MHD equations, a system whose underlying structure is essentially one-dimensional. To the best of our knowledge, no corresponding local exact controllability result has been established for the genuinely multidimensional full compressible MHD system. The presence of the solenoidal magnetic constraint and its interaction with the compressible transport--parabolic structure introduce difficulties that are absent both from the compressible Navier--Stokes equations and from the planar compressible MHD model.

In this paper, we intend to study the local exact controllability for full compressible MHD equations \eqref{1.1} around an equilibrium $(\rb,\ub,\tb,\bb)$ with $\rb,\tb>0$, $\ub\in\mathbb R^d\setminus\{0\}$ and $\bb\in\mathbb R^d$ in a smooth bounded domain $\Omega \subset \mathbb{R}^d$. The velocity, temperature and magnetic controls are their boundary traces on $(0,T)\times\partial\Omega$. 
Let $\bm n$ denote the outward unit normal to $\partial\Omega$. The density control is prescribed only on the incoming part
\[
 \Gamma_-(\bm U)=\{(t,x)\in(0,T)\times\partial\Omega:
                              \bm U(t,x)\cdot\bm n(x)<0\}.
\]
The control time $T>\frac{\ell}{|\ub|}$, where $\ell=\sup_{x\in\Omega}x\cdot e-\inf_{x\in\Omega}x\cdot e$ denotes the directional width of 
$\Omega$ in the direction $e=\ub/|\ub|$.

Our main result can be stated as follows.

\begin{theorem}\label{thm 1.1}
Assume $P$ belongs to $C^3$ locally around $(\bar{\rho}, \bar{\theta})$ and $T>\frac{\ell}{|\ub|}$. 
There exists $\delta>0$ such that for all
$(\rho_0,\bm u_0,\theta_0,\bm b_0)\in H^2(\Omega)\times H^2(\Omega)^d\times H^2(\Omega)\times H^2(\Omega)^d$ satisfying 
\begin{equation}\label{1.2}
 \Div\bm b_0=0,\qquad
 \int_{\Gamma_j}\bm b_0\cdot\bm n\,dS=0
 \quad\text{for every connected component }\Gamma_j\text{ of }\partial\Omega,
\end{equation}
and
$
 \norm{(\rho_0,\bm u_0,\theta_0,\bm b_0)}_{H^2(\Omega)}\le\delta,
$
there exist boundary controls $(v_{\widetilde\rho},\bm v_{\bm u},v_{\widetilde\theta},\bm v_{\bm B})$ and a solution $(\widetilde{\rho},\bm U,\widetilde{\theta},\bm B)$ of \eqref{1.1} satisfying
\[
\begin{aligned}
\bm U=\bm v_{\bm u},\quad \widetilde\theta=v_{\widetilde\theta},\quad
\bm B=\bm v_{\bm B}&&\text{on }(0,T)\times\partial\Omega,\quad\widetilde\rho=v_{\widetilde\rho}&&\text{on }\Gamma_-(\bm U),
\end{aligned}
\]
and the initial condition
\begin{align}\label{1.3}
 (\widetilde\rho,\bm U,\widetilde\theta,\bm B)(0)
   =(\rb+\rho_0,\ub+\bm u_0,\tb+\theta_0,\bb+\bm b_0),  \quad\text{in~}\Omega,
\end{align}
and the final condition
\begin{align}\label{1.4}
 (\widetilde\rho,\bm U,\widetilde\theta,\bm B)(T)
   =(\rb,\ub,\tb,\bb), \quad \text{in~}\Omega.
\end{align}
Moreover, $\widetilde\rho,\widetilde\theta>0$, and the solution $(\widetilde\rho,\bm U,\widetilde\theta,\bm B)$ has the regularity
\begin{align}\label{1.5}
&\widetilde\rho \in C([0,T];H^2(\Omega))\cap H^1(0,T;H^1(\Omega)), \notag\\
&(\bm U,\widetilde\theta,\bm B)\in\big[L^2(0,T;H^3(\Omega))\cap H^1(0,T;H^1(\Omega))\cap C([0,T];H^2(\Omega))\big]^{2d+1}.
\end{align}
The traces of the boundary controls have the regularity 
\begin{align}\label{1.6}
 v_{\widetilde\rho}\in C([0,T];H^{3/2}(\partial\Omega)),~~ (\bm v_{\bm u},v_{\widetilde\theta},\bm v_{\bm B})
 \in\big[L^2(0,T;H^{5/2}(\partial\Omega))
           \cap C([0,T];H^{3/2}(\partial\Omega))\big]^{2d+1}.
\end{align}
\end{theorem}

\begin{remark}
\label{remark 1.1}
The lower bound for the control time $T$ allows the information from waves of density to
travel at velocity (close to) $\ub$, which requires constant states having nonzero velocity. If $\ub=0$,
the systems may not be controllable. One can see \cite{Maity2015} for more studies on the controllability of 1D
linearized compressible Navier-Stokes equations around trivial constant steady state.
\end{remark}

\begin{remark}
\label{remark 1.2}
If $\partial\Omega$ is connected, the flux condition \eqref{1.2} follows from the divergence theorem. For several boundary components it is a necessary compatibility condition for controllability to a constant magnetic field. Indeed, the induction equation and $\Div\bm B=0$ give
$
\frac{d}{dt}\int_{\Gamma_j}\bm B\cdot\bm n\,dS
=\int_{\Gamma_j}\nabla\times(\bm U\times\bm B-\nu\nabla\times\bm B)\cdot\bm n\,dS=0
$
for every closed boundary component $\Gamma_j$, by Stokes' theorem. In two dimensions, the last integral is the integral of a tangential derivative over a closed curve. At the regularity of \eqref{1.5}, the identity follows in distributions in time by the trace theorem and approximation. Since the constant field $\bb$ has zero flux through every $\Gamma_j$, the same must hold for $\bm b_0$. Here, we do not impose zero mean on $\bm b_0$ in $\Omega$.

\end{remark}

\begin{remark}
\label{remark 1.3} 
In contrast with the velocity, temperature, and density controls, which are prescribed as Dirichlet traces (the density only on the inflow set), the magnetic boundary control is not an arbitrary boundary datum. It is the trace of a solenoidal trajectory, such that the divergence constraint 
$\Div\bm B=0$ is preserved. This distinction is specific to the magnetic constraint in the MHD system.
More precisely, an admissible magnetic boundary trace is one which admits a lifting
\[
\bm H\in L^2(0,T; H^3(\Omega)^d)\cap H^1(0,T; H^1(\Omega)^d)\cap C([0,T]; H^2(\Omega)^d)
\]
such that $\Div\bm H=0$ and $\int_{\Gamma_j}\bm H(t)\cdot\bm n\,dS=0$ for every $j$ and $t$. The controls and the solution are constructed together; arbitrary boundary data are not prescribed. In the periodic construction below, $v_{\bm b}$ is the $m$-component magnetic potential control, whereas $\bm v_{\bm B}$ is the $d$-component physical boundary trace.
\end{remark}

The proof of Theorem~\ref{thm 1.1} relies on several structural ingredients that are specific to the multidimensional full compressible MHD system. The first concerns the magnetic constraint $\Div\bm b=0$, which must be preserved. To this end, we construct a solenoidal extension of the magnetic perturbation to the periodic domain and introduce a magnetic potential formulation, such that the localized magnetic control enters through the curl form $\nabla\times(\chi v_{\bm b})$. 
The divergence free condition and the zero spatial mean of the magnetic field were preserved by the flux condition \eqref{1.2} and a Bogovskii-type right inverse.
The second one arises in the controllability analysis of the linearized system. The magnetic coupling links the compressible modes, the solenoidal component of the velocity, and the magnetic variable, and therefore prevents the effective viscous flux reduction from yielding a closed subsystem. We overcome this difficulty by introducing an auxiliary coupled system and combining transport controllability with weighted parabolic estimates based on suitably chosen Carleman weights. The third one is a mixed coordinate formulation of the nonlinear system. While the density, velocity, and temperature are written in transported coordinates to eliminate the nonlinear convection field, the magnetic perturbation is kept in Eulerian coordinates in order to preserve its solenoidal structure. The associated coordinate mismatch terms are shown to be quadratic and can therefore be incorporated into the fixed-point scheme. These constructions together provide a control framework that is compatible with the transport structure, the parabolic dissipation, and the geometric constraint of the magnetic field, and constitute the main analytical ingredients behind the multidimensional controllability result.

The rest of the paper is organized as follows. In section 2, we outline the strategy and methods of the proof of our main results. Section 3 studies the null controllability of the linearized system coupling of parabolic and transport equations. In section 4, by the fixed point argument, we obtain the local exact controllability for MHD flows.

\section{Strategy and methods}\label{section 2}

In order to show the controllability for the full compressible MHD equations in multidimensional case, first of all, let us reduce the boundary control problem to an interior control problem on a torus, as in \cite{EGG2016}. 
Let $ L>0 $ be large enough so that $\overline{\Omega}\subset[0,L]^d$ and no characteristic of the equilibrium flow returns to another periodic copy of $\Omega$ during $[0,T]$, and set $\mathbb{T}_L=(\mathbb{R}/L\mathbb{Z})^d$ and $Q_T=(0,T)\times\mathbb{T}_L$.
To preserve the magnetic constraint, we use a magnetic potential \(a\), which has \(m=1\) component in dimension two and \(m=3\) components in dimension three. Thus \(a\) is \(\mathbb{R}^m\)-valued, with \(m\) determined by the spatial dimension. 
Thus, we can study the control problem in the cube $\mathbb{T}_L$ completed with periodic boundary conditions with the control functions supported in $\mathbb{T}_L\setminus\overline{\Omega}$.

The extended nonlinear system is written in primitive variables $(\widetilde\rho,\bm U,\widetilde\theta,\bm B)$:
\begin{equation}\label{2.1}
\left\{\begin{aligned}
 &\partial_t\widetilde\rho+\Div(\widetilde\rho\bm U)=\check v_\rho, ~~&\text{in~} \mathbb{T}_L,
\\
 &\widetilde\rho(\partial_t\bm U+\bm U\cdot\nabla\bm U)
 -\mu\Delta\bm U-(\lambda'+\mu)\nabla\Div\bm U
 +\nabla P(\widetilde\rho,\widetilde\theta)
 -(\nabla\times\bm B)\times\bm B=\check{\bm v}_{\bm u},~~&\text{in~} \mathbb{T}_L,
\\
 &c_v\widetilde\rho(\partial_t\widetilde\theta+\bm U\cdot\nabla\widetilde\theta)
 -\kappa\Delta\widetilde\theta+P(\widetilde\rho,\widetilde\theta)\Div\bm U
 =2\mu|D(\bm U)|^2+\lambda'(\Div\bm U)^2
   +\nu|\nabla\times\bm B|^2+\check v_\theta,~~&\text{in~}  \mathbb{T}_L,
\\
 &\partial_t\bm B-\nabla\times(\bm U\times\bm B)-\nu\Delta\bm B
 =\nabla\times(\chi v_{\bm b}),\qquad \Div\bm B=0, ~~&\text{in~} \mathbb{T}_L.
\end{aligned}\right.
\end{equation}
The curl control in the last equation is essential for controlling the magnetic perturbation and preserves its divergence and spatial mean.

Secondly, we will use the
following translation
\[
 \check\rho=\widetilde\rho-\rb,\quad \check{\bm u}=\bm U-\ub,
 \quad\check\theta=\widetilde\theta-\tb,\quad\check{\bm b}=\bm B-\bb.
\]
to transform the local exact controllability
of \eqref{2.1} into a null controllability of the system
\begin{equation}\label{2.2}
\left\{\begin{aligned}
 &\big(\partial_t+(\ub+\check{\bm u})\cdot\nabla\big)\check\rho+\rb\Div\check{\bm u}=\check f_\rho+\check v_\rho, ~~&\text{in~} \mathbb{T}_L,\\
 &\rb \big(\partial_t+(\ub+\check{\bm u})\cdot\nabla\big)\check{\bm u}-\mu\Delta\check{\bm u}
 -(\lambda'+\mu)\nabla\Div\check{\bm u}
+\Prb\nabla\check\rho+\Ptb\nabla\check\theta-[(\nabla\times\check{\bm b})\times\bb]
 =\check{\bm f}_{\bm u}+\check{\bm v}_{\bm u}, ~~&\text{in~} \mathbb{T}_L,\\
 &c_v\rb \big(\partial_t+(\ub+\check{\bm u})\cdot\nabla\big)\check\theta-\kappa\Delta\check\theta
 +\Pb\Div\check{\bm u}=\check f_\theta+\check v_\theta, ~~&\text{in~} \mathbb{T}_L,\\
 &(\partial_t+\ub\cdot\nabla)\check{\bm b}-\nu\Delta\check{\bm b}-\nabla\times(\check{\bm u}\times\bb)
=\check{\bm f}_{\bm b}+\nabla\times(\chi v_{\bm b}),\qquad\Div\check{\bm b}=0, ~~&\text{in~} \mathbb{T}_L,
\end{aligned}\right.
\end{equation}
with initial condition
\begin{equation}\label{2.3}
 (\check\rho,\check{\bm u},\check\theta,\check{\bm b})(0)
 =(\check\rho_0,\check{\bm u}_0,\check\theta_0,\check{\bm b}_0),
\end{equation}
where 
$\Pb:=P(\rb,\tb)$, $\Prb:=\partial_\rho P(\rb,\tb)$, and $\Ptb:=\partial_\theta P(\rb,\tb)$ and the nonlinear source terms are exactly
\begin{align*}
 \check f_\rho={}&-\check\rho\Div\check{\bm u},\\
 \check{\bm f}_{\bm u}={}&-\check\rho \big(\partial_t+(\ub+\check{\bm u})\cdot\nabla\big)\check{\bm u}-\nabla\big[P(\rb+\check\rho,\tb+\check\theta)-\Pb
                -\Prb\check\rho-\Ptb\check\theta\big]
 +(\nabla\times\check{\bm b})\times\check{\bm b},\\
 \check f_\theta={}&-c_v\check\rho \big(\partial_t+(\ub+\check{\bm u})\cdot\nabla\big)\check\theta-\big[P(\rb+\check\rho,\tb+\check\theta)-\Pb\big]\Div\check{\bm u}
+2\mu|D(\check{\bm u})|^2+\lambda'(\Div\check{\bm u})^2+\nu|\nabla\times\check{\bm b}|^2,\\
 \check{\bm f}_{\bm b}={}&\nabla\times(\check{\bm u}\times\check{\bm b}).
\end{align*}

The initial perturbations are extended from $\Omega$ to $\T$ such that
\begin{equation}\label{2.4}
 \norm{(\check\rho_0,\check{\bm u}_0,\check\theta_0,\check{\bm b}_0)}_{(H^2(\T))^{2+2d}}
 \le C_L\norm{(\rho_0,\bm u_0,\theta_0,\bm b_0)}_{(H^2(\Omega))^{2+2d}}.
\end{equation}
Here, the density, velocity, and temperature are extended by a standard bounded 
$H^2$ extension. The magnetic field requires a different procedure, because the extension must remain divergence free and have zero spatial mean. This is achieved by the flux condition \eqref{1.2}  and the Bogovskii-type correction described below.

Fix a bounded $H^2$ extension and multiply it by a smooth cutoff equal to one near $\overline\Omega$ and supported in a smooth connected set $K$ with $\overline\Omega\subset\subset K\subset\subset(0,L)^d$, obtaining $\bm e$. Each connected component $D$ of $K\setminus\overline\Omega$ is a bounded smooth domain. Since $g=\Div\bm e\in H^1(K)$ vanishes in $\Omega$ and near $\partial K$, its trace on $\partial D$ is zero, and hence $g|_D\in H^1_0(D)$. It follows from the flux conditions and the divergence theorem that $\int_Dg=0$. Apply \cite[Theorem~4.9]{CM2010} with degree $d$ and $s=1$. On a connected domain, the top-degree compact-support cohomology has dimension one, with the integral as its obstruction. The closed-range assertion therefore gives a bounded linear right inverse of the divergence on $\{g\in H^1_0(D):\int_Dg=0\}$. The supported Sobolev spaces in that theorem agree, by zero extension, with $H^1_0(D)$ and $H^2_0(D)$ on smooth domains. Thus
\[
 \bm z_D\in H^2_0(D)^d,\qquad \Div\bm z_D=g|_D,
 \qquad \norm{\bm z_D}_{H^2(D)}\le C_D\norm g_{H^1(D)}.
\]
Since compactly supported smooth functions are dense in $H^2_0(D)$, the zero extension is bounded on $H^2_0(D)$. Hence $\bm e_\sigma=\bm e-\sum_D\bm z_D$ is a compactly supported $H^2$ solenoidal extension, with no interface contribution to its divergence. Moreover,
$
 \int_K(\bm e_\sigma)_i\,dx
 =\int_K\Div(x_i\bm e_\sigma)\,dx=0.
$
Finally, we periodize the field constructed above, which yields a periodic magnetic field
$\check{\bm b}_0\in H^2_{\sigma,0}(\T)$, where
\[
H^j_{\sigma,0}(\T)
=\{\bm b\in H^j(\T)^d:\Div\bm b=0,\ \int_\T\bm b\,dx=0\}.
\] 
Consequently, \eqref{2.4} holds. Notice that all extension maps are fixed and linear. The flux restriction is automatic only when $\partial\Omega$ is connected and zero mean is imposed on the extension, not on $\bm b_0$ in $\Omega$.

Since the controlled system \eqref{2.2} is still a nonlinear system, in the third step, we need to study the null controllability for a linearized system around $(\bar{\rho}, \bar{u}, \bar{\theta}, \bar{b})$, which reads as follows
\begin{equation}\label{2.5}
\left\{\begin{aligned}
 &(\partial_t+\ub\cdot\nabla)\rho+\rb\Div\bm u=\widehat f_\rho+\chi v_\rho,~~&\text{in~} \mathbb{T}_L,\\
 &\rb(\partial_t+\ub\cdot\nabla)\bm u-\mu\Delta\bm u-(\lambda'+\mu)\nabla\Div\bm u
 +\Prb\nabla\rho+\Ptb\nabla\theta-[(\nabla\times\bm b)\times\bb]
    =\widehat{\bm f}_{\bm u}+\chi\bm v_{\bm u},~~&\text{in~} \mathbb{T}_L,\\
 &c_v\rb(\partial_t+\ub\cdot\nabla)\theta-\kappa\Delta\theta+\Pb\Div\bm u
    =\widehat f_\theta+\chi v_\theta,~~&\text{in~} \mathbb{T}_L,\\
 &(\partial_t+\ub\cdot\nabla)\bm b-\nu\Delta\bm b-\nabla\times(\bm u\times\bb)
    =\widehat{\bm f}_{\bm b}+\nabla\times(\chi v_{\bm b}),\qquad\Div\bm b=0, ~~&\text{in~} \mathbb{T}_L.
\end{aligned}\right.
\end{equation}
Here $\widehat{\bm f}_{\bm b}(t)\in H^1_{\sigma,0}(\T)$, and $\chi\in C^\infty(\T;[0,1])$ is a cutoff function satisfying
\begin{equation}\label{2.6}
 \chi(x)=0\quad\hbox{if }d(x,\Omega)\le\varepsilon,
 \qquad \chi(x)=1\quad\hbox{if }d(x,\Omega)\ge2\varepsilon.
\end{equation}
where $\varepsilon>0$ is small enough, which will be taken in the weights for the Carleman estimates later.

\begin{remark}\label{remark 2.1}
With the help of ${\rm div}\bm b_0=0$, the condition ${\rm div}\bm b=0$ holds naturally
after taking the divergence of this system. Spatial integration separately gives $\partial_t\int_\T\bm b=0$. 
\end{remark}

Taking advantage of the equivalence between the controllability of a linear system and the observability property of its adjoint system, and combining Carleman estimates with suitable weight functions for the parabolic system and the characteristics method for the transport equation, we establish the null controllability of the linearized system \eqref{2.5}. Here two magnetic ingredients are essential and distinguish the present work from the compressible Navier--Stokes case. First, the magnetic control is not imposed as an arbitrary Dirichlet datum. It is generated through a magnetic potential and enters the magnetic equation in the curl form $\nabla\times(\chi v_{\boldsymbol b})$. This is necessary to preserve both the divergence free condition and the zero spatial mean of the magnetic field, as explained in Remark~\ref{remark 2.1}. An arbitrary localized vector source would destroy these constraints and invalidate the functional setting $H^2_{\sigma,0}(\mathbb T_L)$. Second, the magnetic initial data cannot be extended to the torus by the same standard $H^2$ extension used for the density, velocity, and temperature. Since the extension must remain divergence free and zero-mean, we use the flux condition \eqref{1.2} and a Bogovskii-type right inverse to correct the divergence on each connected component of the exterior region. 
If $\partial\Omega$ has several connected components, this flux condition is the compatibility condition explained in Remark~\ref{remark 1.2}.
 These constructions are indispensable for applying the linear controllability theorem on $H^2_{\sigma,0}(\mathbb T_L)$. As for the nonlinear system \eqref{2.2} with initial condition \eqref{2.3}, the quadratic term $\check{u}\cdot\nabla\check{\rho}$ in the density equation causes a loss of regularity, so it cannot be placed in the source term. To establish null controllability, we rewrite the system in new coordinates similar to Lagrangian coordinates for the density, velocity, and temperature so that the nonlinear convective term disappears, while keeping the magnetic field in Eulerian coordinates in order to preserve its divergence free and zero-mean constraints. The resulting composition differences are estimated as quadratic sources, and a Schauder argument closes the nonlinear problem.

\section{Controllability of the linearized system}
\subsection{Preliminaries}
To deal with the magnetic constraint, we introduce the magnetic potential \(a\) instead of working directly with the magnetic field \(\boldsymbol b\) to remain divergence free and zero-mean. Let \(\boldsymbol b=\nabla\times a\). It ensures \(\operatorname{div}\boldsymbol b=0\), and the zero-mean condition is preserved by the curl operator on the torus. Moreover, the magnetic control can be imposed in the curl form \(\nabla\times(\chi v_{\boldsymbol b})\), which remains localized outside \(\Omega\) and preserves both constraints. The following lemma is needed for the construction of the magnetic field.
 \begin{lemma}\label{lemma 3.1}
Let \(j\ge 0\). On the space \(H^j_{\sigma,0}(\T)\) of zero-mean divergence free fields, the following hold.
\begin{enumerate}
\item[(i)] The operator \((-\Delta)^{-1}(\nabla\times)\) is bounded from \(H^j_{\sigma,0}(\T)\) to \(H^{j+1}(\T)^m\), and
\[
(\nabla\times)(-\Delta)^{-1}(\nabla\times)=I
\qquad\text{on } H^j_{\sigma,0}(\T).
\]
\item[(ii)] The operator \((-\Delta)^{-1}(\nabla\times)\) commutes with multiplication by functions of time.
\item[(iii)] Let \(a\) solve a componentwise controlled heat equation with localized source \(\chi_0v\). Then \(\boldsymbol b=\nabla\times a\) satisfies a solenoidal magnetic equation whose source \(\nabla\times(\chi_0v)\) remains localized.
\end{enumerate}
\end{lemma}
\begin{proof}
For $k\in(2\pi/L)\mathbb Z^d\setminus\{0\}$, the Fourier symbols of the curl satisfy
$
(ik\times)(ik\times)=|k|^2I-k\otimes k
$ with $d=3$,
whereas in dimension two the vector-to-scalar and scalar-to-vector symbols are, respectively,
\[
C(k)=i(-k_2,k_1),\qquad K(k)=i(k_2,-k_1)^{\mathsf T},\qquad K(k)C(k)=|k|^2I-k\otimes k.
\]
Since $k\cdot\widehat{\bm b}(k)=0$, both identities imply 
$
(\nabla\times)(-\Delta)^{-1}(\nabla\times)\bm b=\bm b~\text{on }H^j_{\sigma,0}(\T).
$
The lifting symbol is bounded by $|k|^{-1}$; hence
\[
\|(-\Delta)^{-1}(\nabla\times b)\|_{H^{j+1}}\le C\|b\|_{H^j}.
\]
At zero frequency we set the symbol equal to zero and the zero-mean hypothesis makes the identity valid there as well.

Because this operator is a fixed Fourier multiplier, it commutes with multiplication by time dependent functions. The curl operator commutes with both $\partial_t+\bar{\boldsymbol u}\cdot\nabla$ and $\Delta$, and, as a differential operator. Hence, if the source $\chi_0 v$ is localized, then so is $\nabla\times(\chi_0 v)$. The argument uses only this support property of the curl. Finally, applying $\nabla\times$ to the componentwise controlled heat equation for the unconstrained potential $a$, whose source $\chi_0 v$ is localized, gives a solenoidal magnetic equation whose source $\nabla\times(\chi_0 v)$ remains localized.
\end{proof}

For the linear system \eqref{2.5}, write $(\widehat\rho_0,\widehat{\bm u}_0,\widehat\theta_0,\widehat{\bm b}_0)$ for the initial data, with $\widehat{\bm b}_0\in H^2_{\sigma,0}(\T)$. On zero-mean solenoidal fields, let us define
\[
 a_0=(-\Delta)^{-1}(\nabla\times\widehat{\bm b}_0),
 \qquad F_{\bm b}=(-\Delta)^{-1}(\nabla\times\widehat{\bm f}_{\bm b}).
\]
Here $(-\Delta)^{-1}$ is defined by its Fourier multiplier on zero-mean functions. For $\bm b\in H^j_{\sigma,0}$,
\[
 \nabla\times\big[(-\Delta)^{-1}(\nabla\times\bm b)\big]=\bm b,
 \qquad
 \norm{(-\Delta)^{-1}(\nabla\times\bm b)}_{H^{j+1}}
 \le C\norm{\bm b}_{H^j}.
\] 
Thus, we can solve the system obtained from \eqref{2.5} by putting $\bm b=\nabla\times a$ and replacing the last equation by
\begin{equation}\label{3.1}
 (\partial_t+\ub\cdot\nabla) a-\nu\Delta a-(\bm u\times\bb)=F_{\bm b}+\chi v_{\bm b}.
\end{equation}
Applying $\nabla\times$ recovers the magnetic equation in \eqref{2.5}. The condition $a(T)=0$ suffices for $\bm b(T)=0$.

The adjoint of the potential system reads
\begin{equation}\label{3.2}
\left\{\begin{aligned}
 &-(\partial_t+\ub\cdot\nabla)\sigma-\Prb\Div\bm z=g_\sigma, ~~&\text{in~} \mathbb{T}_L,\\
 &-\rb(\partial_t+\ub\cdot\nabla)\bm z-\mu\Delta\bm z-(\lambda'+\mu)\nabla\Div\bm z
       -\rb\nabla\sigma-\Pb\nabla\eta-(\bb\times\phi)=\bm g_{\bm z}, ~~&\text{in~} \mathbb{T}_L,\\
 &-c_v\rb(\partial_t+\ub\cdot\nabla)\eta-\kappa\Delta\eta-\Ptb\Div\bm z=g_\eta, ~~&\text{in~} \mathbb{T}_L,\\
 &-(\partial_t+\ub\cdot\nabla)\phi-\nu\Delta\phi+\nabla\times[\nabla\times(\bm z\times\bb)]=g_\phi, ~~&\text{in~} \mathbb{T}_L.
\end{aligned}\right.
\end{equation}

We know that the effective viscous flux introduced in \cite{Lions1998} is important for the existence theory of compressible flows.
Here, we borrow this kind of variable, i.e.
\begin{equation}\label{3.3}
 q=\Div\bm z+\frac{\rb}{\lambda}\sigma
\end{equation}
to reduce the order of the coupling and derive the observability inequality.
Applying the divergence operator to the second equation of \eqref{3.2} and adding $\rb^2/\lambda$ times its first equation gives
\begin{equation}\label{3.4}
\left\{\begin{aligned}
 &-(\partial_t+\ub\cdot\nabla)\sigma+\frac{\rb\Prb}{\lambda}\sigma-\Prb q=g_\sigma, ~~&\text{in~} \mathbb{T}_L,\\
 &-\rb(\partial_t+\ub\cdot\nabla) q-\lambda\Delta q+\frac{\rb^3\Prb}{\lambda^2}\sigma-\frac{\rb^2\Prb}{\lambda}q -\Pb\Delta\eta-\Div((\bb\times\phi))
       =\Div\bm g_{\bm z}+\frac{\rb^2}{\lambda}g_\sigma, ~~&\text{in~} \mathbb{T}_L,\\
 &-c_v\rb(\partial_t+\ub\cdot\nabla)\eta-\kappa\Delta\eta+\frac{\rb\Ptb}{\lambda}\sigma-\Ptb q=g_\eta,~~&\text{in~} \mathbb{T}_L,\\
 &-\rb(\partial_t+\ub\cdot\nabla)\bm z-\mu\Delta\bm z-(\lambda-\mu)\nabla q
       -\frac{\rb\mu}{\lambda}\nabla\sigma-\Pb\nabla\eta-(\bb\times\phi)=\bm g_{\bm z}, ~~&\text{in~} \mathbb{T}_L,\\
 &-(\partial_t+\ub\cdot\nabla)\phi-\nu\Delta\phi+\nabla\times[\nabla\times(\bm z\times\bb)]=g_\phi, ~~&\text{in~} \mathbb{T}_L.
\end{aligned}\right.
\end{equation}
After we introduce the effective viscous flux \(q\),
the equation for \(q\) only determines the divergence of \(\boldsymbol z\). The solenoidal part of \(\boldsymbol z\), however, is not controlled by \(q\). To retain the full velocity adjoint, we keep \(\boldsymbol z\) as an additional variable and retain its equation, with \(\operatorname{div}\boldsymbol z\) expressed in terms of \(q\) and \(\sigma\).

Thus, one can use duality to obtain the observability property of \eqref{3.4}. Since the system is linear, this property is equivalent to a null controllability statement for the following auxiliary control problem
\begin{equation}\label{3.5}
\left\{\begin{aligned}
 &(\partial_t+\ub\cdot\nabla) r+\frac{\rb\Prb}{\lambda} r=-\frac{\rb^3\Prb}{\lambda^2} y-\frac{\rb\Ptb}{\lambda} h
               -\frac{\rb\mu}{\lambda}\Div\bm w+f_r+\chi_0v_r, ~~&\text{in~} \mathbb{T}_L,\\
 &\rb(\partial_t+\ub\cdot\nabla) y-\lambda\Delta y=\Prb r+\frac{\rb^2\Prb}{\lambda}y+\Ptb h
               -(\lambda-\mu)\Div\bm w+f_y+\chi_0v_y, ~~&\text{in~} \mathbb{T}_L,\\
 &c_v\rb(\partial_t+\ub\cdot\nabla) h-\kappa\Delta h=\Pb\Delta y-\Pb\Div\bm w+f_h+\chi_0v_h, ~~&\text{in~} \mathbb{T}_L,\\
 &\rb(\partial_t+\ub\cdot\nabla)\bm w-\mu\Delta\bm w=-\big\{\bb\times[\nabla\times(\nabla\times\psi)]\big\}+\bm f_{\bm w}+\chi_0\bm v_{\bm w}, ~~&\text{in~} \mathbb{T}_L,\\
 &(\partial_t+\ub\cdot\nabla)\psi-\nu\Delta\psi=-(\nabla y\times\bb)+(\bm w\times\bb)+f_\psi+\chi_0v_\psi, ~~&\text{in~} \mathbb{T}_L.
\end{aligned}\right.
\end{equation}
Here, $\bm w$ is the additional velocity variable and $\chi_0$ is a smooth cutoff function satisfying
\begin{equation}\label{3.6}
\supp\chi_0 \subset\subset\{\chi=1\} \text{~and ~} \chi_0(x)=1 \text{~for all~} x\in \mathbb{T}_L \text{~such that~}
d(x, \Omega)\geq 3\varepsilon.
\end{equation}

Now, let us recall the weight construction used in \cite{BEG2016,EGG2016}. Choose $T_0,T_1,\varepsilon>0$ such that
\begin{equation}\label{3.7}
 T_0+2T_1<T-\frac{\ell+12\varepsilon}{|\ub|}.
\end{equation}
Then, we enlarge the periodicity cell, if necessary, such that a segment of length
$|\ub|T$ meets at most one periodic copy of the $6\varepsilon$-neighborhood
of $\overline\Omega$. The following property follows from \eqref{3.7}.
\begin{lemma}\label{lemma 3.0}
There exist $T_0<t_-<t_+<T-2T_1$ and
$\gamma\in C_c^\infty((T_0,T-2T_1);[0,\infty))$
such that $\gamma\ge1$ on $[t_-,t_+]$ and
\[
M(x):=\int_0^T\gamma(t)\chi_0(x+t\ub)^2\,dt\ge m_0>0
\quad\text{for all }x\in\T.
\]
In particular, $M^{-1}\in C^\infty(\T)$, with bounds depending only on
the fixed geometry and cutoffs.
\end{lemma}
\begin{proof}
In view of \eqref{3.7}, one can take $t_-,t_+$ in the plateau satisfying
$|\ub|(t_+-t_-)>\ell+6\varepsilon$.
For any $x\in\T$, lift the segment $x+t\ub$ to $\mathbb R^d$.
If its distance from the periodic copies of $\Omega$ were less than
$3\varepsilon$ for every $t\in[t_-,t_+]$, the separation of copies would
place the whole segment in a single $3\varepsilon$-neighborhood.
Its displacement in direction $e$ would then be at most
$\ell+6\varepsilon$, a contradiction. Thus
$\chi_0(x+t\ub)=1$ for some $t\in[t_-,t_+]$, by \eqref{3.6}.
By the smoothness of $\gamma$ and continuity, $M(x)>0$ for every
$x\in\T$. From the compactness of $\T$ it follows that
$m_0:=\min_{\T}M>0$. Finally, the assertion about $M^{-1}$ is
obtained by differentiating under the integral sign.
\end{proof}
If necessary, the torus is further enlarged to accommodate the transported weight construction in \cite{EGG2016}.
Let $\widetilde\varphi\in C^2([0,T]\times\T;[6,7])$ satisfy $(\partial_t+\ub\cdot\nabla)\widetilde\varphi=0$ and have no spatial critical points outside a set $\omega \subset\subset\{\chi_0=1\}$, uniformly in time. 
Fix $\lambda_0\ge1$ be fixed sufficiently large and $s\ge s_0\ge1$ be the Carleman parameter, with $s_0\lambda_0^2e^{2\lambda_0}>2$. The geometric quantities above are independent of $s$. Define
\begin{equation}\label{3.8}
 \zeta(t)=\begin{cases}
 1+(1-t/T_0)^{s\lambda_0^2e^{2\lambda_0}},&0\le t\le T_0,\\
 1,&T_0\le t\le T-2T_1,\\
 (T-t)^{-1},&T-T_1\le t<T,
 \end{cases}
\end{equation}
with a smooth increasing connection on $[T-2T_1,T-T_1]$; take $T_1$ small enough that $\zeta\ge1$. Set
\begin{equation}\label{3.9}
 \varphi=\zeta(\lambda_0e^{12\lambda_0}-e^{\lambda_0\widetilde\varphi}),
 \quad \xi=\zeta e^{\lambda_0\widetilde\varphi},
 \quad \Phi=\zeta\lambda_0e^{12\lambda_0}.
\end{equation}

The weights $\zeta$, $\varphi$, $\xi$, and $\Phi$ depend on $s$, which is omitted from the notation. The connection in \eqref{3.8} is chosen independently of $s$. Thus $\zeta(0)=2$, $\frac{14}{15}\Phi\le\varphi\le\Phi$, and $|\zeta'|\le Cs\zeta^2$.
In the contraction estimates, the constant \(C\) is independent of \(s\), whereas the constant \(C_s\) in the higher Sobolev estimates may depend on the value at which \(s\) is eventually fixed.

\subsection{Controllability of some linear equations}
We now begin the controllability analysis. The following standard parabolic regularity estimate will be used throughout this section. In what follows, let us use the abbreviations
\[
L^pH^j=L^p(0,T;H^j(\T)),\qquad H^1H^j=H^1(0,T;H^j(\T)),\qquad CH^j=C([0,T];H^j(\T)).
\]
$\norm{\cdot}_2$ denotes the $L^2(Q_T)$ norm. For a trace at a fixed time, it denotes the $L^2(\T)$ norm.
For the Riesz representations below, we use the Fourier norms
\[
\norm v_{H^q(\T)}^2=\sum_{k\in(2\pi/L)\mathbb Z^d}(1+|k|^2)^q|\widehat v(k)|^2,
\]
where the Fourier coefficients are taken in the orthonormal basis of $L^2(\T)$.
 
For \(c,d_0>0\), the parabolic equation
\[
c(\partial_t+\bar{\boldsymbol u}\cdot\nabla)v-d_0\Delta v=F,\qquad v(0)=v_0,
\]
satisfies, for every \(r\in\mathbb R\),
\begin{equation}\label{3.10}
\|v\|_{L^2H^{r+2}}+\|\partial_t v\|_{L^2H^r}+\|v\|_{C([0,T];H^{r+1})}
\le C\bigl(\|F\|_{L^2H^r}+\|v_0\|_{H^{r+1}}\bigr).
\end{equation}
This is a standard parabolic regularity estimate. Removing the drift by the Fourier multiplier \(e^{i(\bar{\boldsymbol u}\cdot k)t}\), the equation reduces to a family of scalar ODEs in \(k\). Young's inequality and the heat kernel decay give the \(H^{r+2}\) bound, the zero mode is treated separately, and summation with Sobolev weights yields the spatial estimate. The time derivative follows from the equation, and the trace bound follows by interpolation. The same argument applies to the Lamé system, since its symbol
\[
\mu|k|^2I+(\lambda'+\mu)k\otimes k
\]
is positive definite for $k\neq0$ by the assumptions $\mu>0$ and $\lambda=2\mu+\lambda'>0$, and diagonalizes on \(k^\perp\) and \(\mathbb R k\). 
This regularity estimate will be combined with the Carleman estimate below to obtain the weighted controllability estimates for the scalar heat equation.

\begin{lemma}\label{lemma 3.2}
Let $a_1>0$. For sufficiently large $\lambda_0$ and $s$, every smooth $W$ satisfies
\begin{align}\label{3.11}
 &s^{3/2}\lambda_0^2\norm{\xi^{3/2}We^{-s\varphi}}_2
 +s^{1/2}\lambda_0\norm{\xi^{1/2}\nabla We^{-s\varphi}}_2
 +s\lambda_0^{3/2}e^{7\lambda_0}\norm{W(0)e^{-s\varphi(0)}}_2\notag\\
 &\qquad\le C\norm{(-(\partial_t+\ub\cdot\nabla)-a_1\Delta)We^{-s\varphi}}_2
       +Cs^{3/2}\lambda_0^2\norm{\xi^{3/2}\chi_0We^{-s\varphi}}_2.
\end{align}
The same weights can be used for any fixed finite collection of positive diffusion coefficients.
\end{lemma}
\begin{proof}
The estimate is the scalar Carleman estimate established in \cite{BEG2016,EGG2016}. To handle the drift term, we write
$
-\partial_t-a_1\Delta=\bigl(-(\partial_t+\ub\cdot\nabla)-a_1\Delta\bigr)+\ub\cdot\nabla.
$
Applying the drift-free estimate to the first operator and using \(\xi\ge e^{6\lambda_0}\), we find that the drift term is absorbed by the gradient term for sufficiently large \(s\). Taking the maximum of the finitely many thresholds yields the last assertion. 
\end{proof}

We now apply the Carleman estimate to the scalar parabolic equation.
\begin{lemma}\label{lemma 3.3}
Fix $c_0,d_0>0$. There is a linear choice of $(Y,v_Y)$ solving
\begin{equation}\label{3.12}
 c_0(\partial_t+\ub\cdot\nabla) Y-d_0\Delta Y=F+\chi_0v_Y,
 \qquad Y(0)=Y_0,\qquad Y(T)=0,
\end{equation}
such that, for $Y_0\in H^1$ and $\norm{\zeta^{-3/2}Fe^{s\varphi}}_2<\infty$,
\begin{align}\label{3.13}
 &s^{3/2}\norm{Ye^{s\varphi}}_2
 +s^{1/2}\norm{\zeta^{-1}\nabla Ye^{s\varphi}}_2
 +s^{-1/2}\norm{\zeta^{-2}\nabla^2Ye^{s\varphi}}_2\notag\\
 &\quad+\norm{\zeta^{-3/2}\chi_0v_Ye^{s\varphi}}_2
 \le C\big(\norm{\zeta^{-3/2}Fe^{s\varphi}}_2
       +s^{1/2}\norm{Y_0e^{s\varphi(0)}}_2
       +s^{-1/2}\norm{\nabla Y_0e^{s\varphi(0)}}_2\big).
\end{align}
For $Y_0\in L^2$, the same assertion holds without the Hessian and the initial gradient terms. In addition, this choice satisfies
\begin{equation}\label{3.14}
 \norm{\chi_0v_Ye^{\frac67s\Phi}}_{L^2H^2}
 \le C_s\big(\norm{\zeta^{-3/2}Fe^{s\varphi}}_2+\norm{Y_0e^{s\Phi(0)}}_{L^2}\big),
\end{equation}
and, for $j=1,2$,
\begin{align}\label{3.15}
 \norm{Ye^{\frac67s\Phi}}_{L^2H^{j+2}}
 \le C_s\big(&\norm{Fe^{\frac67s\Phi}}_{L^2H^j}
       +\norm{\zeta^{-3/2}Fe^{s\varphi}}_2+\norm{Y_0e^{s\Phi(0)}}_{H^{j+1}}\big).
\end{align}
\end{lemma}
\begin{proof}
Set
$
L=c_0(\partial_t+\ub\cdot\nabla)-d_0\Delta$ and $L^*=-c_0(\partial_t+\ub\cdot\nabla)-d_0\Delta.
$
On smooth test functions, define the bilinear form
\[
B(W,W)=\norm{e^{-s\varphi}L^*W}_2^2
+s^3\norm{\chi_0\zeta^{3/2}e^{-s\varphi}W}_2^2.
\]
By Lemma~\ref{lemma 3.2}, the linear functional
\[
\ell_Y(W)
= \int_{Q_T} F W \,dx\,dt
+ c_0 \int_{\mathbb{T}_L} Y_0 W(0) \,dx
\]
satisfies
\[
|\ell_Y(W)|
\le C\bigl(s^{-3/2}\norm{\zeta^{-3/2}Fe^{s\varphi}}_2
+s^{-1}\norm{Y_0e^{s\varphi(0)}}_2\bigr)B(W,W)^{1/2}.
\]
The bilinear form $B$ is positive definite on smooth test functions. Indeed, $B(W,W)=0$ implies $W=0$ by Lemma~\ref{lemma 3.2}. Completing this space with respect to $B$ yields a Hilbert space $\mathcal H$. Since $\ell_Y$ is bounded on $\mathcal H$, the Riesz representation theorem provides a unique $W_*\in\mathcal H$ such that $B(W_*,V)=\ell_Y(V)$ for all $V\in\mathcal H$, together with
\begin{equation}\label{3.16}
B(W_*,W_*)^{1/2}
\le C\bigl(s^{-3/2}\norm{\zeta^{-3/2}Fe^{s\varphi}}_2
+s^{-1}\norm{Y_0e^{s\varphi(0)}}_2\bigr).
\end{equation}
The weighted residual $e^{-s\varphi}L^*W_*$ and observation $\chi_0\zeta^{3/2}e^{-s\varphi}W_*$ are understood as the $L^2$ limits of those of a smooth approximating sequence. We then define
\[
Y=e^{-2s\varphi}L^*W_*
,\qquad
v_Y=-s^3\zeta^3\chi_0e^{-2s\varphi}W_*.
\]
This choice depends linearly on \((Y_0,F)\). Moreover, the Riesz identity yields \eqref{3.12}, including both endpoint conditions, in transposition. Since \(0\le\chi_0\le1\), estimate \eqref{3.16} further gives
\begin{equation}\label{3.17}
s^{3/2}\norm{Ye^{s\varphi}}_2
+\norm{\zeta^{-3/2}\chi_0v_Ye^{s\varphi}}_2
\le C\bigl(\norm{\zeta^{-3/2}Fe^{s\varphi}}_2
+s^{1/2}\norm{Y_0e^{s\varphi(0)}}_2\bigr).
\end{equation}

We now turn to the weighted gradient estimate. Let
\[
q=e^{s\varphi}Y,\qquad
H=e^{s\varphi}(F+\chi_0v_Y).
\]
A direct differentiation gives
\begin{equation}\label{3.18}
(c_0(\partial_t+\ub\cdot\nabla)-d_0\Delta)q
=H-2d_0s\nabla\varphi\cdot\nabla q
+\bigl(c_0s(\partial_t+\ub\cdot\nabla)\varphi
+d_0s^2|\nabla\varphi|^2-d_0s\Delta\varphi\bigr)q.
\end{equation}
With \(\lambda_0\) fixed, we have
\[
|\nabla\varphi|+|\Delta\varphi|\le C\zeta,\qquad
|(\partial_t+\ub\cdot\nabla)\varphi|\le Cs\zeta^2,\qquad
|\zeta'|\le Cs\zeta^2.
\]
Denote the first- and zero-order coefficients by \(B_1\) and \(B_0\), respectively, such that
$
|B_1|\le Cs\zeta, |B_0|\le Cs^2\zeta^2.
$
Multiplying \eqref{3.18} by \(s\zeta^{-2}q\) and integrating, we observe that the spatial drift \(\ub\) contributes nothing, since its coefficient and \(\zeta\) are independent of the spatial variable. Young's inequality then yields
\begin{align*}
&s\left|\int\zeta^{-2}Hq\right|
\le C\norm{\zeta^{-3/2}H}_2^2+Cs^2\norm q_2^2, ~~\left|\int\zeta^{-2}B_1\cdot\nabla q\,q\right|
\le\tfrac{d_0s}{2}\norm{\zeta^{-1}\nabla q}_2^2+Cs^3\norm q_2^2,\\
&s\left|\int\zeta^{-2}B_0q^2\right|
+s\left|\int(\zeta^{-2})'q^2\right|
\le Cs^3\norm q_2^2,
\end{align*}
where we used \(\zeta\ge1\) and \(s\ge1\). Consequently,
\[
s\norm{\zeta^{-1}\nabla q}_2^2
\le C\bigl(\norm{\zeta^{-3/2}H}_2^2+s^3\norm q_2^2+s\norm{q(0)}_2^2\bigr).
\]
Since \(e^{s\varphi}\nabla Y=\nabla q-sq\nabla\varphi\), estimate \eqref{3.17} implies the required gradient bound in \eqref{3.13}.

To obtain the weighted Hessian bound, assume \(Y_0\in H^1\) and multiply \eqref{3.18} by \(-s^{-1}\zeta^{-4}\Delta q\). Integrating by parts in time for the gradient energy, we have
\begin{align*}
s^{-1}\norm{\zeta^{-2}\Delta q}_2^2
\le C\bigl(&s^{-1}\norm{\zeta^{-2}H}_2^2
+s\norm{\zeta^{-1}\nabla q}_2^2+s^3\norm q_2^2 +\norm{\zeta^{-3/2}\nabla q}_2^2
+s^{-1}\norm{\nabla q(0)}_2^2\bigr).
\end{align*}
To verify the powers, note that
$
s^{-1}|B_1|^2\zeta^{-4}\le Cs\zeta^{-2}, 
s^{-1}|B_0|^2\zeta^{-4}\le Cs^3,
$
while differentiating \(\zeta^{-4}\) produces a bound \(Cs\zeta^{-3}\) before multiplication by \(s^{-1}\). On the torus, \(\norm{D^2q}_2\le C\norm{\Delta q}_2\). Furthermore,
\[
e^{s\varphi}D^2Y
=D^2q-s(\nabla\varphi\otimes\nabla q+\nabla q\otimes\nabla\varphi)
+(s^2\nabla\varphi\otimes\nabla\varphi-sD^2\varphi)q.
\]
Combining the zeroth- and first-order estimates therefore yields the asserted \(s^{-1/2}\zeta^{-2}\) Hessian bound. At \(t=0\), we have
\[
\nabla q(0)=e^{s\varphi(0)}\nabla Y_0+s\nabla\varphi(0)e^{s\varphi(0)}Y_0,
\]
so the extra contribution is bounded precisely by
\[
s^{1/2}\norm{Y_0e^{s\varphi(0)}}_2
+s^{-1/2}\norm{\nabla Y_0e^{s\varphi(0)}}_2.
\]
These energy calculations are justified first with temporal cutoffs and smooth data; the bounds, lower semicontinuity, and the transposition identity then pass to the completion. They do not presume that the exponentially rescaled state has a classical terminal value.

It remains to prove the regularity of the control. Set
$
V=e^{-\frac{211}{210}s\Phi}W_*.
$
In the distributional sense,
\[
(-c_0(\partial_t+\ub\cdot\nabla)-d_0\Delta)V
=e^{-\frac{211}{210}s\Phi}(-c_0(\partial_t+\ub\cdot\nabla)-d_0\Delta)W_*
+c_0\frac{211}{210}s\Phi'e^{-\frac{211}{210}s\Phi}W_*.
\]
The first term equals \(e^{s\varphi-\frac{211}{210}s\Phi}(Ye^{s\varphi})\) and is bounded in \(L^2\). The second is controlled by the Carleman estimate on \(W_*\) and the strict inequality \(\frac{211}{210}>1\). A smooth test-field approximation gives \(V(T)=0\), and the backward version of \eqref{3.10} yields
\[
\norm V_{L^2H^2\cap H^1L^2}
\le C_s\bigl(\norm{\zeta^{-3/2}Fe^{s\varphi}}_2
+\norm{Y_0e^{s\varphi(0)}}_2\bigr).
\]
Since
\[
e^{\frac67s\Phi}\chi_0v_Y
=-s^3\zeta^3\chi_0^2
e^{(\frac67+\frac{211}{210})s\Phi-2s\varphi}V,
\]
and the coefficient together with its first two spatial derivatives is bounded by \(C_s\), multiplication is bounded on \(H^2\). This proves \eqref{3.14} for the same \(v_Y\).

Finally, applying the parabolic regularity estimate \eqref{3.10} with \(r=j\) to the equation for
$
Y_*=e^{\frac67s\Phi}Y
$
proves \eqref{3.15}, with constants depending on the fixed value of \(s\).
\end{proof}

The following lemma extends Lemma~\ref{lemma 3.3} to weights with powers of \(\zeta\). The substitution \(Z=\zeta^k Y\) reduces the proof to Lemma~\ref{lemma 3.3}, since the extra term \(c_0 k(\zeta'/\zeta)Z\) is absorbed for large \(s\); the estimates transfer directly, and dividing by \(\zeta^k\) gives the stated bounds. The case \(k=1\) is needed for the magnetic coupling.
\begin{lemma}\label{lemma 3.4}
For $k=1/2$ and $k=1$, the same equation and endpoint conditions as in Lemma~\ref{lemma 3.3} admit a linear choice satisfying
\begin{align*}
 &s^{3/2}\norm{\zeta^kYe^{s\varphi}}_2
 +s^{1/2}\norm{\zeta^{k-1}\nabla Ye^{s\varphi}}_2
 +s^{-1/2}\norm{\zeta^{k-2}\nabla^2Ye^{s\varphi}}_2 +\norm{\zeta^{k-3/2}\chi_0v_Ye^{s\varphi}}_2\\
 &\le C\big(\norm{\zeta^{k-3/2}Fe^{s\varphi}}_2
 +s^{1/2}\norm{Y_0e^{s\varphi(0)}}_2
 +s^{-1/2}\norm{\nabla Y_0e^{s\varphi(0)}}_2\big).
\end{align*}
For $L^2$ initial data, omit the Hessian and initial-gradient terms. The regularity estimates are explicitly
\begin{align*}
 \norm{\chi_0v_Ye^{\frac67s\Phi}}_{L^2H^2}
 &\le C_s\big(\norm{\zeta^{k-3/2}Fe^{s\varphi}}_2
                    +\norm{Y_0e^{s\Phi(0)}}_{L^2}\big),\\
 \norm{Ye^{\frac67s\Phi}}_{L^2H^{j+2}}
 &\le C_s\big(\norm{Fe^{\frac67s\Phi}}_{L^2H^j}
       +\norm{\zeta^{k-3/2}Fe^{s\varphi}}_2 +\norm{Y_0e^{s\Phi(0)}}_{H^{j+1}}\big),\quad j=1,2.
\end{align*}
\end{lemma}

With the weight function in Carleman estimate at hand, one can study the control problem for the transport equation.
\begin{lemma}\label{lemma 3.5}
Under \eqref{3.7}, there is a linear choice of $(r,v_r)$ satisfying
\begin{equation}\label{3.19}
 (\partial_t+\ub\cdot\nabla) r+\frac{\rb\Prb}{\lambda} r=F+\chi_0v_r,
 \qquad r(0)=r_0,\qquad r(T)=0,
\end{equation}
with
\begin{equation}\label{3.20}
 \norm{\zeta^{-1}(r,\chi_0v_r)e^{s\varphi}}_2
 \le C\big(\norm{\zeta^{-1}Fe^{s\varphi}}_2
                    +\norm{r_0e^{s\varphi(0)}}_2\big).
\end{equation}
For $j=1,2$ this selection also satisfies
\begin{equation}\label{3.21}
 \norm{(r,\chi_0v_r)e^{\frac67s\Phi}}_{L^2H^j}
 \le C_s\big(\norm{Fe^{\frac67s\Phi}}_{L^2H^j}
                  +\norm{r_0e^{s\Phi(0)}}_{H^j}\big).
\end{equation}
\end{lemma}
\begin{proof}
Set \(c=\frac{\rb\Prb}{\lambda}\). By Lemma~\ref{lemma 3.0}, there exists a compact interval \([t_-,t_+]\) contained in the plateau such that every characteristic meets \(\{\chi_0>0\}\) within this interval.
 Choose a nonnegative smooth function \(\gamma\), supported in the plateau and positive on \([t_-,t_+]\), and define
\[
M(x)=\int_0^T\gamma(t)\chi_0(x+t\ub)^2\,dt>0.
\]
By compactness, \(M\ge m_0>0\) for some constant \(m_0\), and the derivatives of \(M^{-1}\) are bounded. We then set
\begin{align*}
&B(x)=r_0(x)+\int_0^T e^{c\tau}F(\tau,x+\tau\ub)\,d\tau, ~~
h(t,x)=M(x)^{-1}\int_0^t\gamma(\tau)\chi_0(x+\tau\ub)^2\,d\tau,\\
&v_r(t,x+t\ub)=-e^{-ct}\gamma(t)\chi_0(x+t\ub)M(x)^{-1}B(x),\\
&r(t,x+t\ub)=e^{-ct}\left[r_0(x)+\int_0^t e^{c\tau}F(\tau,x+\tau\ub)\,d\tau-h(t,x)B(x)\right].
\end{align*}
These directly yield \eqref{3.19}, its endpoint conditions, and the linear dependence on the data.

Next, let
$
A(x)=\lambda_0e^{12\lambda_0}-e^{\lambda_0\widetilde\varphi(0,x)}>0,
$
such that \(\varphi(t,x+t\ub)=\zeta(t)A(x)\). For sufficiently large \(s\), the function \(z^{-1}e^{szA(x)}\) is increasing on \(z\ge1\), uniformly in \(x\). Before the plateau, we have \(h=0\) and \(\zeta\) is decreasing, whence
\[
\frac{\zeta(\tau)}{\zeta(t)}
e^{s\varphi(t,x+t\ub)-s\varphi(\tau,x+\tau\ub)}\le1,
\qquad 0\le\tau\le t,
\]
on that interval. After the control support, where \(h=1\), we rewrite the state as
\[
r(t,x+t\ub)=-e^{-ct}\int_t^T e^{c\tau}F(\tau,x+\tau\ub)\,d\tau.
\]
On the final increasing interval, the same quotient is at most one for \(t\le\tau\). On the plateau, the weight equals \(e^{sA(x)}\), which is its minimum along the characteristic. Consequently, the weighted forward and backward kernels are bounded by \(e^{|c|T}\), with \(L^2\) operator norms at most \(Te^{|c|T}\). Moreover,
\[
e^{sA(x)}|B(x)|
\le e^{s\varphi(0,x)}|r_0(x)|
+e^{|c|T}\int_0^T\zeta(\tau)^{-1}e^{s\varphi(\tau,x+\tau\ub)}
|F(\tau,x+\tau\ub)|\,d\tau.
\]
Using the boundedness of \(h\), \(\gamma\chi_0/M\), and \(e^{\pm ct}\), and integrating in \(x\) and \(t\), we obtain \eqref{3.20} uniformly in large \(s\); note that translations preserve the spatial measure.

Finally, for \eqref{3.21}, we observe that \(e^{\frac67s\Phi(t)}\) has the same monotonicity on the three intervals. Differentiating the characteristic formulas up to order \(j\), we see that the spatial derivatives of \(h\) and of the control coefficient are bounded and supported on the plateau. Hence the same kernel estimates apply to every Leibniz term, which proves \eqref{3.21}. Smooth approximation preserves the formulas and the transport traces.
\end{proof}

\subsection{Controllability of the auxiliary system \eqref{3.5}}
Taking advantage of the null controllability for the linear parabolic equation and the transport equation, we can obtain the controllability for
the auxiliary system \eqref{3.5}. In this subsection, the spatial variable lies in \(\T\).
\begin{theorem}\label{thm 3.6}
Let \eqref{3.7} hold. For sufficiently large fixed $s$, initial data
\[
 r_0\in H^2,\qquad (y_0,h_0,\bm w_0,\psi_0)
              \in H^3\times H^3\times(H^3)^d\times(H^3)^m,
\]
and sources $\bm f=(f_r,f_y,f_h,\bm f_{\bm w},f_\psi)$ with
$e^{s\Phi}\bm f\in L^2(H^2\times H^2\times H^2\times(H^2)^d\times(H^2)^m)$, there is a linear selection of states and controls for \eqref{3.5} such that
\begin{equation}\label{3.22}
 (r,y,h,\bm w,\psi)(T)=0.
\end{equation}
It satisfies
\begin{align}\label{3.23}
 &\norm{re^{\frac67s\Phi}}_{L^2H^2}
 +\norm{(y,h,\bm w,\psi)e^{\frac67s\Phi}}_{L^2(H^4\times H^4\times(H^4)^d\times(H^4)^m)}\notag\\
 &\quad+\norm{\chi_0(v_r,v_y,v_h,\bm v_{\bm w},v_\psi)e^{\frac67s\Phi}}_{L^2(H^2\times H^2\times H^2\times(H^2)^d\times(H^2)^m)}\notag\\
 &\le C_s\big(\norm{e^{s\Phi}\bm f}_{L^2(H^2\times H^2\times H^2\times(H^2)^d\times(H^2)^m)}
 +\norm{r_0e^{s\Phi(0)}}_{H^2}\notag\\
 &\hspace{18mm}+\norm{(y_0,h_0,\bm w_0,\psi_0)e^{s\Phi(0)}}_{H^3\times H^3\times(H^3)^d\times(H^3)^m}\big).
\end{align}
\end{theorem}
\begin{proof}
We want to use the controllability for the linear equations
and the Banach fixed point theorem to prove this theorem. Let us first introduce 
a work space $\mathcal Z_s$ consisting of periodic distributional tuples $(r,y,h,\bm w,\psi)$ of component dimensions $(1,1,1,d,m)$, for which the following norm is finite, with all spatial derivatives understood in the weak sense:
\begin{equation}\label{3.24}
\begin{aligned}
\norm{(r,y,h,\bm w,\psi)}_{\mathcal Z_s}=& s^{1/8}\norm{\zeta^{-1}re^{s\varphi}}_2+\big(s^{3/2}\norm{\zeta^{1/2}ye^{s\varphi}}_2
       +s^{1/2}\norm{\zeta^{-1/2}\nabla ye^{s\varphi}}_2
       +s^{-1/2}\norm{\zeta^{-3/2}\nabla^2ye^{s\varphi}}_2\big)\\
 &+\big(s^{3/4}\norm{he^{s\varphi}}_2
       +s^{-1/4}\norm{\zeta^{-1}\nabla he^{s\varphi}}_2
       +s^{-5/4}\norm{\zeta^{-2}\nabla^2he^{s\varphi}}_2\big)\\
 &+\big(s^{5/4}\norm{\zeta^{1/2}\bm we^{s\varphi}}_2
       +s^{1/4}\norm{\zeta^{-1/2}\nabla\bm we^{s\varphi}}_2
       +s^{-3/4}\norm{\zeta^{-3/2}\nabla^2\bm we^{s\varphi}}_2\big)\\
 &+\big(s^{15/8}\norm{\zeta\psi e^{s\varphi}}_2
       +s^{7/8}\norm{\nabla\psi e^{s\varphi}}_2
       +s^{-1/8}\norm{\zeta^{-1}\nabla^2\psi e^{s\varphi}}_2\big).
\end{aligned}
\end{equation}
The equations are solved diagonally with all coupling terms evaluated at a given input in $\mathcal Z_s$.
 Use Lemma~\ref{lemma 3.5} for $r$, Lemma~\ref{lemma 3.4} with $k=1/2$ for $y$ and $\bm w$, Lemma~\ref{lemma 3.3} for $h$, and Lemma~\ref{lemma 3.4} with $k=1$ for $\psi$. 
The resulting map has fixed initial data and zero terminal data. Denote its state component by $\mathcal T_s$. The scalar source norms, in the order $(r,y,h,\bm w,\psi)$, carry the factors
$
s^{1/8}\zeta^{-1}, \zeta^{-1}, 
s^{-3/4}\zeta^{-3/2}, s^{-1/4}\zeta^{-1}, s^{3/8}\zeta^{-1/2},
$
respectively, in addition to $e^{s\varphi}$. Substitution of the coupling terms in \eqref{3.5} into these norms gives \eqref{3.25}, and shows that $\mathcal T_s$ maps $\mathcal Z_s$ into itself.

We record every coupling bound. For the difference of two inputs, denote the five nonnegative grouped terms in \eqref{3.24} by $R,Y,H,W,\Psi$, respectively, and denote the corresponding output terms by primes. Lemmas~\ref{lemma 3.3}, \ref{lemma 3.4}, and \ref{lemma 3.5}, together with $\zeta\ge1$, give
\begin{align}\label{3.25}
 R'&\le C(s^{-11/8}Y+s^{-5/8}H+s^{-1/8}W), ~~~
 Y'\le C(s^{-1/8}R+s^{-3/2}Y+s^{-3/4}H+s^{-1/4}W),\notag\\
 H'&\le C(s^{-1/4}Y+s^{-1}W), ~~~
 W'\le Cs^{-1/8}\Psi, ~~~
 \Psi'\le C(s^{-1/8}Y+s^{-7/8}W). 
\end{align}
For the two second-order couplings and the transport derivative, we have
\begin{align*}
 s^{-3/4}\norm{\zeta^{-3/2}\Delta y\,e^{s\varphi}}_2
 &\le Cs^{-1/4}Y, ~~~
 s^{-1/4}\norm{\zeta^{-1}\big\{\bb\times[\nabla\times(\nabla\times\psi)]\big\}\,e^{s\varphi}}_2
 \le Cs^{-1/8}\Psi,\\
 s^{1/8}\norm{\zeta^{-1}\Div\bm w\,e^{s\varphi}}_2
 &\le Cs^{-1/8}W,
\end{align*}
and the first-order source in the last equation is estimated by
\[
 s^{3/8}\norm{\zeta^{-1/2}(\nabla y\times\bb)\,e^{s\varphi}}_2
 \le Cs^{-1/8}Y.
\]
The remaining zero- and first-order terms account for the other powers displayed in \eqref{3.25}. Consequently, the Lipschitz constant of the diagonal map is at most \(Cs^{-1/8}\). 
Choosing $s$ sufficiently large that $Cs^{-1/8}\le1/2$, and above all the scalar thresholds, therefore guarantees a unique fixed point in $\mathcal Z_s$ for this choice of scalar operators. Since these selections are fixed and linear, uniqueness of the fixed point preserves linear dependence on the data.
The initial and inhomogeneous terms are controlled by the right-hand side of \eqref{3.23}, and the fixed point together with its controls depends linearly on the data.

It remains to establish \eqref{3.23}, since the Sobolev orders are required in the subsequent duality argument. To begin with, the basic norm yields \(y,h,\bm w,\psi\in L^2H^2\) after multiplication by any time weight \(e^{\alpha s\Phi}\) with \(\alpha<14/15\). It also bounds the four parabolic controls in \(e^{\frac67s\Phi}L^2H^2\) by Lemmas~\ref{lemma 3.3} and~\ref{lemma 3.4}. Since the last equation of \eqref{3.5} has an \(H^1\) source, \eqref{3.15} first gives \(\psi e^{\frac67s\Phi}\in L^2H^3\). The equation for \(\bm w\) then implies \(\bm w e^{\frac67s\Phi}\in L^2H^3\). As a consequence, the source of the transport equation belongs to \(e^{-\frac67s\Phi}L^2H^2\), and Lemma~\ref{lemma 3.5} yields \(r e^{\frac67s\Phi}\) and \(\chi_0v_r e^{\frac67s\Phi}\) in \(L^2H^2\).

The equation for \(y\) now has an \(H^2\) source and hence gives \(y e^{\frac67s\Phi}\in L^2H^4\). In succession, the equations for \(h\), \(\psi\), and \(\bm w\) provide the same \(H^4\) regularity. All these implications are estimates whose right-hand sides are bounded by the preceding step and by the data norms appearing on the right of \eqref{3.23}. Derivatives of the time weights are handled using the strict gap \(6/7<14/15\), which is already present in the basic norm. This completes the proof of \eqref{3.23}.
\end{proof}

\subsection {Controllability of the linear system \eqref{2.5}}
In this part, we will use the relationship between observability and controllability to derive the following controllability result for system \eqref{2.5}.
\begin{theorem}\label{thm 3.7}
Assume \eqref{3.7}. Let
$
 (\widehat\rho_0,\widehat{\bm u}_0,\widehat\theta_0)\in H^2\times(H^2)^d\times H^2,
~ \widehat{\bm b}_0\in H^2_{\sigma,0}(\T).
$
For all sufficiently large fixed $s$, suppose
$
 \widehat f_\rho e^{s\Phi}\in L^2H^2, ~
 (\widehat{\bm f}_{\bm u},\widehat f_\theta)e^{\frac76s\Phi}
                 \in L^2((H^1)^d\times H^1), ~
 \widehat{\bm f}_{\bm b}e^{s\Phi}\in L^2H^1_{\sigma,0}.
$
Then \eqref{2.5} admits controls $(v_\rho,\bm v_{\bm u},v_\theta,v_{\bm b})$ and a solution depending linearly on these data, such that
\begin{equation}\label{3.26}
 (\rho,\bm u,\theta,\bm b)(T)=0.
\end{equation}
Moreover, if we set
\begin{align*}
 \mathcal D_s={}&\norm{\widehat f_\rho e^{s\Phi}}_{L^2H^2}
 +\norm{(\widehat{\bm f}_{\bm u},\widehat f_\theta)e^{\frac76s\Phi}}_{L^2((H^1)^d\times H^1)}
 +\norm{\widehat{\bm f}_{\bm b}e^{s\Phi}}_{L^2H^1}\\
 &+e^{\frac76s\Phi(0)}\norm{(\widehat\rho_0,\widehat{\bm u}_0,\widehat\theta_0,\widehat{\bm b}_0)}_{H^2\times(H^2)^d\times H^2\times(H^2)^d},
\end{align*}
then we have
\begin{align}\label{3.27}
 \norm{(\rho,\bm u,\theta,\bm b)e^{\frac67s\Phi}}_{L^2(H^2\times(H^3)^d\times H^3\times(H^3)^d)} +\norm{\chi(v_\rho,\bm v_{\bm u},v_\theta,v_{\bm b})e^{\frac67s\Phi}}_{L^2(H^2\times(H^1)^d\times H^1\times(H^2)^m)}
 \le C_s\mathcal D_s.
\end{align}
The magnetic source actually applied in \eqref{2.5} is $\nabla\times(\chi v_{\bm b})$, and it is supported outside $\overline\Omega$. 
\end{theorem}
\begin{proof}
Let \(p=(\sigma,\bm z,\eta,\phi)\) be an arbitrary smooth test field, and let
\(g=(g_\sigma,\bm g_{\bm z},g_\eta,g_\phi)\) denote the four residuals obtained by
substituting \(p\) into the left-hand sides of \eqref{3.2}. Define \(q\) by \eqref{3.3}
and set \(Z=(\sigma,q,\eta,\bm z,\phi)\). The residual in \eqref{3.4} is then given by
$
G=\left(g_\sigma,\Div\bm g_{\bm z}+\frac{\rb^2}{\lambda}g_\sigma,
g_\eta,\bm g_{\bm z},g_\phi\right).
$
For notational brevity, we set
\[
\mathbb H_j=H^j\times H^j\times H^j\times(H^j)^d\times(H^j)^m,
\quad \mathbb J=H^2\times H^3\times H^3\times(H^3)^d\times(H^3)^m,
\]
\[
\mathbb K=H^2\times H^4\times H^4\times(H^4)^d\times(H^4)^m.
\]
Let \(U=(r,y,h,\bm w,\psi)\), \(v=(v_r,v_y,v_h,\bm v_{\bm w},v_\psi)\), and
\(M_A=\operatorname{diag}(1,\rb,c_v\rb,\rb I_d,I_m)\). Integrating by parts and
using \eqref{3.22} yields
\[
\int_{Q_T}\bm f\cdot Z+\dual{M_AU(0)}{Z(0)}
=\int_{Q_T}U\cdot G-\int_{Q_T}\chi_0v\cdot Z.
\]
Taking the supremum over sources in \(e^{-s\Phi}L^2\mathbb H_2\) and initial data in
\(\mathbb J\), it follows from Theorem~\ref{thm 3.6} that
\begin{align}\label{3.28}
\norm{Ze^{-s\Phi}}_{L^2\mathbb H_{-2}}+\norm{Z(0)}_{\mathbb J'} \le C_s\left(
\norm{\chi_0Ze^{-\frac67s\Phi}}_{L^2\mathbb H_{-2}}
+\norm{Ge^{-\frac67s\Phi}}_{L^2\mathbb K'}\right).
\end{align}
The constant masses \((1,\rb,c_v\rb,\rb,1)\) in the initial pairing yield equivalent
norms, which are absorbed into \(C_s\).

We now introduce the observation and residual norms
\begin{align}\label{3.29}
O_s(p)^2={}&\norm{\chi\sigma e^{-\frac67s\Phi}}_{L^2H^{-2}}^2
+\norm{\chi\bm z e^{-\frac67s\Phi}}_{L^2H^{-1}}^2
+\norm{\chi\eta e^{-\frac67s\Phi}}_{L^2H^{-1}}^2
+\norm{\chi\phi e^{-\frac67s\Phi}}_{L^2H^{-2}}^2,\notag\\
N_s(g)^2={}&\norm{g_\sigma e^{-\frac67s\Phi}}_{L^2H^{-2}}^2
+\norm{\bm g_{\bm z}e^{-\frac67s\Phi}}_{L^2H^{-3}}^2
+\norm{g_\eta e^{-\frac67s\Phi}}_{L^2H^{-3}}^2
+\norm{g_\phi e^{-\frac67s\Phi}}_{L^2H^{-4}}^2.
\end{align}
Since \(\chi=1\) in a neighborhood of \(\supp\chi_0\), the distributional identity
\(\chi_0\Div\bm z=\Div(\chi_0\chi\bm z)-\nabla\chi_0\cdot\chi\bm z\),
together with \eqref{3.3}, yields
\begin{equation}\label{3.30}
\norm{\chi_0q}_{H^{-2}}
\le C\big(\norm{\chi\bm z}_{H^{-1}}+\norm{\chi\sigma}_{H^{-2}}\big).
\end{equation}
We observe that the residual \(G\) lies in \(\mathbb K'\). Consequently, the right-hand side of \eqref{3.28} is bounded by \(C_s(O_s(p)+N_s(g))\).

An additional derivative is required for the velocity and temperature adjoints. To this end, we consider \(\bm z e^{-\frac76s\Phi}\) and \(\eta e^{-\frac76s\Phi}\), both of which possess vanishing terminal values. Their backward heat equations, derived from \eqref{3.4}, have right-hand sides belonging to \(L^2H^{-3}\). For \(\bm z e^{-\frac76s\Phi}\), these right-hand sides are given by
\[
\left((\lambda-\mu)\nabla q+\frac{\rb\mu}{\lambda}\nabla\sigma
+\Pb\nabla\eta+(\bb\times\phi)+\bm g_{\bm z}\right)e^{-\frac76s\Phi},
\]
together with the term generated by differentiating the time weight. Estimate \eqref{3.28} controls the latter term since the function \(|s\Phi'|e^{-s\Phi/6}\) is bounded for fixed \(s\). Similarly, the equation for \(\eta e^{-\frac76s\Phi}\) has as its source
$
\left(\Ptb q-\frac{\rb\Ptb}{\lambda}\sigma+g_\eta\right)e^{-\frac76s\Phi},
$
along with the corresponding time derivative term. Consequently, maximal regularity for the backward heat equation yields
\begin{equation}\label{3.31}
\norm{(\bm z,\eta)e^{-\frac76s\Phi}}_{L^2H^{-1}}
+\norm{(\bm z(0),\eta(0))e^{-\frac76s\Phi(0)}}_{H^{-2}}
\le C_s\big(O_s(p)+N_s(g)\big).
\end{equation}
Here \eqref{3.10} with \(r=-3\) is applied separately to each of the two positive diffusion operators. Combining \eqref{3.28} and \eqref{3.31}, we obtain
\begin{align}\label{3.32}
&\norm{\sigma e^{-s\Phi}}_{L^2H^{-2}}
+\norm{(\bm z,\eta)e^{-\frac76s\Phi}}_{L^2H^{-1}}
+\norm{\phi e^{-s\Phi}}_{L^2H^{-2}}\notag\\
&\quad+\norm{\sigma(0)}_{H^{-2}}
+\norm{(\bm z(0),\eta(0))}_{H^{-2}}
+\norm{\phi(0)}_{H^{-3}}
\le C_s\big(O_s(p)+N_s(g)\big).
\end{align}
The initial exponential factors may either be written explicitly or absorbed into \(C_s\).

We conclude the proof with a weighted Hilbert-space construction. By
Lemma~\ref{lemma 3.1}, we have
\[
\norm{a_0}_{H^3}\le C\norm{\widehat{\bm b}_0}_{H^2},\qquad
\norm{F_{\bm b}e^{s\Phi}}_{L^2H^2}
\le C\norm{\widehat{\bm f}_{\bm b}e^{s\Phi}}_{L^2H^1}.
\]
Define the data functional
\begin{align*}
\ell(p)={}\int_{Q_T}\big(\widehat f_\rho\sigma
+\widehat{\bm f}_{\bm u}\cdot\bm z+\widehat f_\theta\eta+F_{\bm b}\cdot\phi\big)+\int_\T\big(\widehat\rho_0\sigma(0)
+\rb\widehat{\bm u}_0\cdot\bm z(0)
+c_v\rb\widehat\theta_0\eta(0)+a_0\cdot\phi(0)\big).
\end{align*}
Here, the potential pairing is Euclidean, including the case \(m=1\). By \eqref{3.32},
\[
|\ell(p)|\le C_s\mathcal D_s(O_s(p)+N_s(g)).
\]
Now quotient smooth test fields by the kernel of
\[
\norm p_{\mathcal P_s}^2=O_s(p)^2+N_s(g)^2.
\]
Minimizing \(\frac12\|p\|_{P_s}^2-\ell(p)\) over this quotient yields a unique minimizer \(p^{*}\); denote its residual by \(g^{*}\). The minimizer satisfies \(\|p^{*}\|_{P_s}\le C_s\mathcal D_s\).

To pass from the minimizer to state, set \(U=(\rho,\mathbf u,\theta,a)\), \(v=(v_\rho,\mathbf v_{\mathbf u},v_\theta,v_{\mathbf b})\), \(\mathbf r=(2,3,3,4)\) and \(\mathbf o=(2,1,1,2)\). Here $r_i$ denotes an entry of this fixed index vector, not the transport unknown $r$ in \eqref{3.5}. Since
\((I-\Delta)^{-j}:H^{-j}\to H^j\) is the Riesz isometry for the Fourier norms fixed above,
 the state and control
representatives are, componentwise,
\begin{align}
U_i &= e^{-\frac{12}{7}s\Phi}(I-\Delta)^{-r_i}g_i^*, \qquad i=1,\ldots,4, \label{3.33} \\
v_i &= -e^{-\frac{12}{7}s\Phi}(I-\Delta)^{-o_i}(\chi p_i^*), \qquad i=1,\ldots,4. \label{3.34}
\end{align}
Their \(e^{\frac67s\Phi}L^2H^{r_i}\) and \(e^{\frac67s\Phi}L^2H^{o_i}\) norms are
exactly the corresponding residual and observation norms, hence bounded by
\(C_s\mathcal D_s\). The applied control is \(\chi v\), supported outside
\(\overline\Omega\); no locality of \((I-\Delta)^{-j}\) is assumed.

It remains to interpret the above construction as a weak solution of the potential system. Since the residual and observation maps extend continuously to \(\mathcal P_s\), with values in the weighted negative Sobolev spaces that define \(N_s\) and \(O_s\), their values at \(p^{*}\) are well defined by completion. For every smooth \(p\), the Euler equation reads
\[
\int_{Q_T}U\cdot g=\ell(p)+\int_{Q_T}\chi v\cdot p,
\]
which is precisely the potential system in transposition. From it, one obtains \((\partial_t+\ub\cdot\nabla)\rho\in L^2H^2\) and
\(\partial_t(\bm u,\theta,a)\in L^2(H^1\times H^1\times H^2)\). Transport regularity
and \eqref{3.10} then provide traces in \(H^2\) for the density and in
\(H^2,H^2,H^3\) for the other entries. Integrating by parts with independently
prescribed test values at \(0\) and \(T\) identifies these traces with the initial
data and zero. Setting \(\mathbf b=\nabla\times a\) then produces both magnetic traces and gives \eqref{3.26}, \eqref{3.27}, as well as the support and constraint properties, by Lemma 3.1. Finally, the Riesz construction determines a single selection, and this selection depends linearly on all the data.
\end{proof}

\begin{remark}\label{remark 3.1}
Theorem~\ref{thm 3.7} also gives the time regularity needed when restricting controlled trajectories to $\Omega$:
\begin{align}\label{3.35}
 &\norm{\rho e^{\frac56s\Phi}}_{C([0,T];H^2)}
 +\norm{\partial_t(\rho e^{\frac56s\Phi})}_{L^2H^1}+\norm{(\bm u,\theta,\bm b)e^{\frac56s\Phi}}_{L^2H^3\cap H^1H^1\cap C([0,T];H^2)}
 \le C_s\mathcal D_s.
\end{align}
As in \cite{EGG2016}, estimate \eqref{3.35} follows by applying the transport and parabolic regularity estimates to the equations for $(\rho,\bm u,\theta,\bm b)e^{\frac56s\Phi}$.
\end{remark}

\section {Proof of Theorem~\ref{thm 1.1}} 
We prove the main result by combining a coordinate transformation, the linear controllability result of Section 3, and Schauder's fixed point theorem. The density, velocity and temperature are written in transported coordinates, whereas the magnetic perturbation remains in Eulerian coordinates. This preserves the solenoidal, zero-mean condition on the magnetic source. Fixing $s$ as in Theorem~\ref{thm 3.7}, we begin by introducing a fixed-point map. 
\subsection{Definition of the fixed point map}
For $W=(\rho,\bm u,\theta,\bm b)$, define
\begin{align}\label{4.1}
 \norm W_{\mathcal B_s}=\norm{e^{5s\Phi/6}W}_{L^\infty H^2}
 +\norm{e^{5s\Phi/6}(\bm u,\theta,\bm b)}_{L^2H^3}+\norm{e^{5s\Phi/6}(\partial_tW+\ub\cdot\nabla W)}_{L^2H^1},
\end{align}
and set
\begin{equation}\label{4.2}
 \mathcal C_R=\{W:\norm W_{\mathcal B_s}\le R,\quad W(T)=0,
              \quad\Div\bm b=0,\quad\textstyle\int_\T\bm b\,dx=0\}.
\end{equation}
Let $\mathcal G_{\rm st}$ denote the state component of the linear selection from Theorem~\ref{thm 3.7}. Then
\begin{equation}\label{4.3}
 \norm{\mathcal G_{\rm st}(W_0,f)}_{\mathcal B_s}
 \le C_s\bigl(\norm{W_0}_{H^2}+\norm f_{\mathcal Y_s}\bigr),
\end{equation}
where $\mathcal Y_s$ is the space of sources with finite norm and magnetic constraint
\[
 \norm f_{\mathcal Y_s}
 =\norm{e^{s\Phi}f_\rho}_{L^2H^2}
 +\norm{e^{7s\Phi/6}(\bm f_{\bm u},f_\theta)}_{L^2H^1}
 +\norm{e^{s\Phi}\bm f_{\bm b}}_{L^2H^1},\qquad
 \bm f_{\bm b}(t)\in H^1_{\sigma,0}.
\]
The first two state norms follow from \eqref{3.35}. The equations in \eqref{2.5} give the last norm, because their right-hand sides have the required $e^{5s\Phi/6}L^2H^1$ bounds by \eqref{3.27}.  

For $W\in\mathcal C_R$, define
\begin{equation}\label{4.4}
 Y(t,x)=x-\int_t^T\bm u(\tau,x+(\tau-t)\ub)\,d\tau.
\end{equation}
It satisfies
$
 \partial_tY+\ub\cdot\nabla Y=\ub+\bm u, Y(T,x)=x.
$
Let $Z(t,\cdot)$ denote its inverse. Thus $Y(t,Z(t,x))=x$ and $Z(t,Y(t,x))=x$. The ordinary derivatives of the inverse satisfy
\begin{equation}\label{4.5}
 \sum_{j=1}^d(\partial_iZ_j)(t,Y(t,x))\partial_jY_k(t,x)=\delta_{ik},
 \qquad (D_xZ)(t,Y(t,x))=(D_xY(t,x))^{-1}.
\end{equation}
Here $(D_xY)_{ij}=\partial_jY_i$. The symbol $I$ denotes the identity map or identity matrix, as appropriate; $Y-I$ and $Z-I$ are periodic displacements.

\begin{proposition}\label{prop 4.1}
There exists $R_0>0$ such that, for $W\in\mathcal C_R$ and $R\le R_0$, $Y,Z$ are uniformly bi-Lipschitz $H^3$ diffeomorphisms of the torus and
\begin{align}\label{4.6}
 \norm{e^{5s\Phi/12}(Y-I)}_{CH^3}
 +\norm{e^{5s\Phi/12}(Z-I)}_{CH^3} +\norm{e^{5s\Phi/12}\bigl[(D_xZ)(t,Y(t,\cdot))-I\bigr]}_{CH^2}
 \le C_sR.
\end{align}
Moreover $\chi(Z(t,x))=0$ for $x\in\overline\Omega$, and $Y(T,x)=Z(T,x)=x$.
\end{proposition}
\begin{proof}
Differentiating \eqref{4.4}, one can get the transport equation. Since translations preserve Sobolev norms and \(\Phi(\tau)\ge\Phi(t)/2\) for \(\tau\ge t\), we obtain
\[
e^{5s\Phi(t)/12}\norm{Y(t)-I}_{H^3}
\le\int_t^T e^{5s\Phi(\tau)/6}\norm{\bm u(\tau)}_{H^3}\,d\tau
\le\sqrt{T-t}\norm{e^{5s\Phi/6}\bm u}_{L^2(t,T;H^3)}.
\]
The Bochner integral is continuous into \(H^3\) away from \(T\), and this bound provides its weighted continuous extension by zero at \(T\). Shrinking \(R_0\) if necessary, we may assume \(\norm{D_xY-I}_\infty\le1/4\). On the periodic lift, the map \(x\mapsto y-(Y(t,x)-x)\) is a contraction; its fixed point therefore defines an inverse \(Z\) on the torus, with uniform bi-Lipschitz and Jacobian bounds.

Next, the \(H^2\) algebra property and smooth matrix inversion near \(I\) give
\[
\norm{(D_xY)^{-1}-I}_{H^2}\le C\norm{Y-I}_{H^3}.
\]
Differentiating the identity \(Y(t,Z(t,x))=x\) then yields
\begin{align*}
D_xZ(t,x)&=[(D_xY)(t,Z(t,x))]^{-1},\\
D_x^2Z(t,x)&=-[(D_xY)(t,Z(t,x))]^{-1}
(D_x^2Y)(t,Z(t,x))[D_xZ(t,x),D_xZ(t,x)].
\end{align*}
By change of variables, \(D_x^2Z\) is bounded in \(L^2\cap L^6\). Differentiating once more produces a term containing \((D_x^3Y)(t,Z(t,x))\) together with products of \((D_x^2Y)(t,Z(t,x))\) and \(D_x^2Z\). The latter lie in \(L^2\) by virtue of the embedding \(L^6\cdot L^3\subset L^2\). Combining this with the identity \(Z(t,x)-x=-[Y(t,Z(t,x))-Z(t,x)]\), we deduce
\[
\norm{Z-I}_{H^3}\le C\bigl(\norm{Y-I}_{H^3}+\norm{Y-I}_{H^3}^2\bigr)\le C\norm{Y-I}_{H^3}.
\]
Smooth approximation justifies the identities and the continuity of inversion for \(H^3\) maps. Multiplication by \(e^{5s\Phi(t)/12}\) then proves \eqref{4.6}, including the terminal limit.

Finally, we decrease \(R_0\) further so that \(\norm{Z-I}_\infty<\varepsilon/2\). If \(\operatorname{dist}(x,\Omega)<\varepsilon/4\), then \(\operatorname{dist}(Z(t,x),\Omega)<3\varepsilon/4\), and \eqref{2.6} implies \(\chi(Z(t,x))=0\). Consequently, the reconstructed controls vanish near \(\overline\Omega\), and the terminal identities follow by inversion.
\end{proof}

The physical perturbations are recovered by
\begin{equation}\label{4.7}
 \check\rho(t,x)=\rho(t,Z(t,x)), ~
 \check{\bm u}(t,x)=\bm u(t,Z(t,x)), ~ \check\theta(t,x)=\theta(t,Z(t,x)),~
 \check{\bm b}(t,x)=\bm b(t,x).
\end{equation}
For any one of the first three components $f(t,x)=\check f(t,Y(t,x))$, the chain rule reads
\begin{align*}
 \partial_tf+\ub\cdot\nabla f
 &=(\partial_t\check f)(t,Y(t,x))
   +(\nabla\check f)(t,Y(t,x))\cdot
         (\partial_tY+\ub\cdot\nabla Y)(t,x)\\
 &=\bigl[\partial_t\check f+(\ub+\check{\bm u})\cdot\nabla\check f\bigr](t,Y(t,x)),\\
 (\partial_i\check f)(t,Y(t,x))
 &=\sum_{j=1}^d(\partial_iZ_j)(t,Y(t,x))\partial_jf(t,x),\\
 (\partial_{il}\check f)(t,Y(t,x))
 &=\sum_{j,k=1}^d(\partial_iZ_j)(t,Y(t,x))(\partial_lZ_k)(t,Y(t,x))\partial_{jk}f(t,x) +\sum_{k=1}^d(\partial_{il}Z_k)(t,Y(t,x))\partial_kf(t,x).
\end{align*}
This removes the variable transport field from the density equation. The magnetic field is not changed to these coordinates. For simplicity, throughout the rest of this section all summations are from \(1\) to \(d\).
A term such as $(\partial_iZ_j)(t,Y(t,x))$ always means: differentiate $Z_j$ with respect to its spatial argument first, and then evaluate at $Y(t,x)$.

Thus, applying the above identities to \eqref{2.2} and separating the linear
part, we reduce our problem to the controllability of the following
system
\begin{equation}\label{4.9}
\left\{\begin{aligned}
 &\partial_t\rho+\ub\cdot\nabla\rho+\rb\Div\bm u=f_\rho(W)+\chi v_\rho,\\
 &\rb(\partial_t\bm u+\ub\cdot\nabla\bm u)-\mu\Delta\bm u
  -(\lambda'+\mu)\nabla\Div\bm u+\Prb\nabla\rho+\Ptb\nabla\theta\\
 &\hspace{15mm}-(\nabla\times\bm b)\times\bb=\bm f_{\bm u}(W)+\chi\bm v_{\bm u},\\
 &c_v\rb(\partial_t\theta+\ub\cdot\nabla\theta)-\kappa\Delta\theta
                +\Pb\Div\bm u=f_\theta(W)+\chi v_\theta,\\
 &\partial_t\bm b+\ub\cdot\nabla\bm b-\nu\Delta\bm b
             -\nabla\times(\bm u\times\bb)
      =\bm f_{\bm b}(W)+\nabla\times(\chi v_{\bm b}),\qquad\Div\bm b=0.
\end{aligned}\right.
\end{equation}
Here the first three equations are written in reference coordinates,
whereas the fourth remains in Eulerian coordinates, with both spatial
variables denoted by $x$. The first three controls satisfy
\[
(\check v_\rho,\check{\mathbf v}_u,\check v_\theta)(t,Y(t,x))
=
\chi(x)(v_\rho,\mathbf v_u,v_\theta)(t,x).
\]
The magnetic control remains $\nabla\times(\chi v_b)$ in Eulerian
coordinates. The initial and terminal conditions are
$W(0)=\mathcal I(W)$, with $\mathcal I(W)$ defined in \eqref{4.14}, and
$W(T)=0$, respectively. The source terms are as follows:
\begin{align}\label{4.10}
 f_\rho(W)={}&-\rho\sum_{i,j}(\partial_iZ_j)(t,Y(t,x))\partial_ju_i-\rb\sum_{i,j}\bigl[(\partial_iZ_j)(t,Y(t,x))-\delta_{ij}\bigr]\partial_ju_i.
\end{align}
\begin{align}\label{4.11}
 f_{\bm u,i}(W)={}&-\rho(\partial_tu_i+\ub\cdot\nabla u_i)+\mu\sum_{j,k}\bigl[(\nabla Z_j\cdot\nabla Z_k)(t,Y(t,x))-\delta_{jk}\bigr]\partial_{jk}u_i\notag\\
 &+\mu\sum_k(\Delta Z_k)(t,Y(t,x))\partial_ku_i +(\lambda'+\mu)\sum_{l,j,k}\bigl[(\partial_iZ_j)(t,Y(t,x))
                (\partial_lZ_k)(t,Y(t,x))-\delta_{ij}\delta_{lk}\bigr]\partial_{jk}u_l\notag\\
 &+(\lambda'+\mu)\sum_{l,k}(\partial_{il}Z_k)(t,Y(t,x))\partial_ku_l -\sum_j(\partial_iZ_j)(t,Y(t,x))
       \partial_j\bigl[P(\rb+\rho,\tb+\theta)-\Pb-\Prb\rho-\Ptb\theta\bigr]\notag\\
 &-\Prb\sum_j\bigl[(\partial_iZ_j)(t,Y(t,x))-\delta_{ij}\bigr]\partial_j\rho -\Ptb\sum_j\bigl[(\partial_iZ_j)(t,Y(t,x))-\delta_{ij}\bigr]\partial_j\theta\notag\\
 &+\bigl[(\nabla\times\bm b)(t,Y(t,x))\times\bb-(\nabla\times\bm b)(t,x)\times\bb\bigr]_i +\bigl[(\nabla\times\bm b)(t,Y(t,x))\times\bm b(t,Y(t,x))\bigr]_i.
\end{align}
\begin{align}\label{4.12}
 f_\theta(W)={}&-c_v\rho(\partial_t\theta+\ub\cdot\nabla\theta)+ \kappa\sum_{j,k}\bigl[(\nabla Z_j\cdot\nabla Z_k)(t,Y(t,x))-\delta_{jk}\bigr]\partial_{jk}\theta +\kappa\sum_k(\Delta Z_k)(t,Y(t,x))\partial_k\theta\notag\\
 &-\Pb\sum_{i,j}\bigl[(\partial_iZ_j)(t,Y(t,x))-\delta_{ij}\bigr]\partial_ju_i -\bigl[P(\rb+\rho,\tb+\theta)-\Pb\bigr]
            \sum_{i,j}(\partial_iZ_j)(t,Y(t,x))\partial_ju_i\notag\\
 &+\frac\mu2\sum_{i,j}\left[\sum_k\bigl((\partial_iZ_k)(t,Y(t,x))\partial_ku_j
                    +(\partial_jZ_k)(t,Y(t,x))\partial_ku_i\bigr)\right]^2\notag\\
 &+\lambda'\left[\sum_{i,j}(\partial_iZ_j)(t,Y(t,x))\partial_ju_i\right]^2
                   +\nu|(\nabla\times\bm b)(t,Y(t,x))|^2.
\end{align}
\begin{equation}\label{4.13}
 \bm f_{\bm b}(W)=\nabla\times\Bigl[
   \bigl(\bm u(t,Z(t,x))-\bm u(t,x)\bigr)\times\bb
                  +\bm u(t,Z(t,x))\times\bm b(t,x)\Bigr].
\end{equation}

Let $(\check\rho_0,\check{\bm u}_0,\check\theta_0,\check{\bm b}_0)$ be the extensions of Section~\ref{section 2}. Set
\begin{equation}\label{4.14}
 \mathcal I(W)(x)=\bigl(\check\rho_0(Y(0,x)),\check{\bm u}_0(Y(0,x)),
                    \check\theta_0(Y(0,x)),\check{\bm b}_0(x)\bigr).
\end{equation}
The magnetic component stays in $H^2_{\sigma,0}$. Define
\begin{equation}\label{4.15}
 \mathcal F(W)=\mathcal G_{\rm st}\bigl(\mathcal I(W),f(W)\bigr).
\end{equation}
Proposition~\ref{prop 4.1} provides the coordinate transformation. Lemma~\ref{lemma 4.3} below will verify the source assumptions of Theorem~\ref{thm 3.7}, and hence show that \eqref{4.15} is well-defined.

\subsection{Estimates for the nonlinear sources}
To estimate the nonlinear source terms, we shall use the following basic Sobolev product and composition estimates.
\begin{lemma}\label{lemma 4.2}
For \(d=2,3\), it holds that
\begin{equation}\label{4.16}
\norm{fg}_{H^2}\le C\norm f_{H^2}\norm g_{H^2},\qquad
\norm{fg}_{H^1}\le C\norm f_{H^2}\norm g_{H^1}.
\end{equation}
If \(V=I+v\) is an \(H^3\) diffeomorphism sufficiently close to the identity in \(H^3\), then \(g\mapsto g(V(\cdot))\) is uniformly bounded on \(H^j\), \(j=1,2,3\), and
\begin{equation}\label{4.17}
\norm{g(V(\cdot))-g}_{H^j}\le C\norm v_{H^3}\norm g_{H^{j+1}},
\qquad j=1,2.
\end{equation}
For the maps in Proposition~\ref{prop 4.1}, the required coordinate differences satisfy
\begin{align}\label{4.18}
&\left\|\sum_{i,j}\bigl[(\partial_iZ_j)(t,Y(t,\cdot))-\delta_{ij}\bigr]\partial_ju_i\right\|_{H^2}
\le C\norm{Y-I}_{H^3}\norm{\bm u}_{H^3},\notag\\
&\left\|\sum_j\bigl[(\partial_iZ_j)(t,Y(t,\cdot))-\delta_{ij}\bigr]\partial_jq\right\|_{H^1}
\le C\norm{Y-I}_{H^3}\norm q_{H^2},\notag\\
&\left\|\sum_{j,k}\bigl[(\nabla Z_j\cdot\nabla Z_k)(t,Y(t,\cdot))-\delta_{jk}\bigr]\partial_{jk}g
+\sum_k(\Delta Z_k)(t,Y(t,\cdot))\partial_kg\right\|_{H^1}\le C\norm{Y-I}_{H^3}\norm g_{H^3},\\
&\left\|\sum_{l,j,k}\bigl[(\partial_iZ_j)(t,Y(t,\cdot))(\partial_lZ_k)(t,Y(t,\cdot))
-\delta_{ij}\delta_{lk}\bigr]\partial_{jk}u_l+\sum_{l,k}(\partial_{il}Z_k)(t,Y(t,\cdot))\partial_ku_l\right\|_{H^1}
\le C\norm{Y-I}_{H^3}\norm{\bm u}_{H^3}.\notag
\end{align}
The estimates with free index \(i\) hold for every \(i=1,\ldots,d\).
\end{lemma}
\begin{proof}
It follows from the embeddings in two and three dimensions that \eqref{4.16} holds.
The chain rule, the uniform Jacobian bounds, and Sobolev multiplication give
$\norm{g(V(\cdot))}_{H^j}\le C\norm g_{H^j}$ for $j=1,2,3$, with $C$ uniform for $V-I$ in a sufficiently small $H^3$ ball. The assertion follows for general $g\in H^j$ by density. For smooth $g$, the fundamental theorem of calculus gives
\[
 g(x+v(x))-g(x)=\int_0^1(Dg)(x+\tau v(x))\cdot v(x)\,d\tau.
\]
Every map $x\mapsto x+\tau v(x)$ has uniform Jacobian bounds. For $j=1$, multiply the $H^1$ function $(Dg)(\cdot+\tau v(\cdot))$ by $v\in H^2$. For $j=2$, use the $H^2$ algebra inequality. Integration in $\tau$ proves \eqref{4.17}. Note that the extra derivative on $g$ is essential. We shall not apply its $H^2$ estimate to arbitrary fixed $H^2$ initial data.

Turning to the coordinate differences,
one can find
\begin{align*}
 &\sum_{i,j}(\partial_iZ_j)(t,Y(t,x))\partial_ju_i-\Div\bm u
     =\sum_{i,j}\bigl[(\partial_iZ_j)(t,Y(t,x))-\delta_{ij}\bigr]\partial_ju_i,\\
 &\sum_j(\partial_iZ_j)(t,Y(t,x))\partial_jq-\partial_iq
     =\sum_j\bigl[(\partial_iZ_j)(t,Y(t,x))-\delta_{ij}\bigr]\partial_jq.
\end{align*}
It follows from \(Z(t,Y(t,x))=x\) that
$
(D_x Z)(t,Y(t,x)) = (D_x Y)^{-1}.
$
By Proposition~\ref{prop 4.1}, one has
\[
\|(D_x Z)(t,Y(t,x)) - I\|_{H^2} \le C \|Y - I\|_{H^3}.
\]
Multiplication in $H^2$ gives the divergence estimate, and $H^2\cdot H^1\to H^1$ gives the gradient estimate.

Taking $l=i$ and summing over $i$ in the chain-rule identity following \eqref{4.7}, we obtain
\begin{align*}
 (\Delta\check g)(t,Y(t,x))-\Delta g(t,x)
  ={}&\sum_{j,k}\bigl[(\nabla Z_j\cdot\nabla Z_k)(t,Y(t,x))-\delta_{jk}\bigr]\partial_{jk}g(t,x)\\
   &+\sum_k(\Delta Z_k)(t,Y(t,x))\partial_kg(t,x),
\end{align*}
where $g(t,x)=\check g(t,Y(t,x))$ and
\begin{align*}
 &(\nabla Z_j\cdot\nabla Z_k)(t,Y(t,x))-\delta_{jk} =\sum_i\Bigl\{
  \delta_{ij}\bigl[(\partial_iZ_k)(t,Y(t,x))-\delta_{ik}\bigr]
 +\delta_{ik}\bigl[(\partial_iZ_j)(t,Y(t,x))-\delta_{ij}\bigr]\\
 &\hspace{50mm}
 +\bigl[(\partial_iZ_j)(t,Y(t,x))-\delta_{ij}\bigr]
  \bigl[(\partial_iZ_k)(t,Y(t,x))-\delta_{ik}\bigr]\Bigr\}.
\end{align*}
By Proposition~\ref{prop 4.1}, its $H^2$ norm is bounded by $C\norm{Y-I}_{H^3}$. For the lower-order coefficient, using the inverse-map estimate, one has
\[
 \norm{(\Delta Z_k)(t,Y(t,\cdot))}_{H^1}
 \le C\norm{Z-I}_{H^3}\le C\norm{Y-I}_{H^3}.
\]
Now $D^2g\in H^1$ and $Dg\in H^2$. Multiplying the principal coefficient by $D^2g$, and the lower-order coefficient by $Dg$, proves the third line of \eqref{4.18}. 

Similarly, for the gradient of the divergence, we have
\begin{align*}
 &(\partial_i\Div\check{\bm u})(t,Y(t,x))-\partial_i\Div\bm u(t,x)\\
 &=\sum_{l,j,k}\bigl[(\partial_iZ_j)(t,Y(t,x))(\partial_lZ_k)(t,Y(t,x))
          -\delta_{ij}\delta_{lk}\bigr]\partial_{jk}u_l(t,x) +\sum_{l,k}(\partial_{il}Z_k)(t,Y(t,x))\partial_ku_l(t,x),
\end{align*}
where
\begin{align*}
(\partial_iZ_j)(t,Y(t,x))(\partial_lZ_k)(t,Y(t,x))-\delta_{ij}\delta_{lk} = &\delta_{ij}\bigl[(\partial_lZ_k)(t,Y(t,x))-\delta_{lk}\bigr]
       +\delta_{lk}\bigl[(\partial_iZ_j)(t,Y(t,x))-\delta_{ij}\bigr]\\
 & +\bigl[(\partial_iZ_j)(t,Y(t,x))-\delta_{ij}\bigr]
              \bigl[(\partial_lZ_k)(t,Y(t,x))-\delta_{lk}\bigr].
\end{align*}
Its $H^2$ norm is bounded by $C\norm{Y-I}_{H^3}$ and we also have
\[
 \norm{(\partial_{il}Z_k)(t,Y(t,\cdot))}_{H^1}
 \le C\norm{Z-I}_{H^3}\le C\norm{Y-I}_{H^3}.
\]
Apply $H^2\cdot H^1\to H^1$ to the principal coefficient times $D^2\bm u$, and to $D\bm u$ times the lower-order coefficient. This gives $C\norm{Y-I}_{H^3}\norm{\bm u}_{H^3}$, completing \eqref{4.18}.
\end{proof}

Now, we are ready to estimate the nonhomogeneous terms.
\begin{lemma}\label{lemma 4.3}
For $R\le R_0$ and $W\in\mathcal C_R$,
\begin{align}\label{4.19}
 \norm{e^{s\Phi}f_\rho(W)}_{L^2H^2}
 +\norm{e^{7s\Phi/6}(\bm f_{\bm u}(W),f_\theta(W))}_{L^2H^1}
 +\norm{e^{s\Phi}\bm f_{\bm b}(W)}_{L^2H^1}&\le C_sR^2,\notag\\
 \bm f_{\bm b}(W)(t)&\in H^1_{\sigma,0}(\T)\quad\text{for a.e. }t\in(0,T).
\end{align}
Moreover, it holds that $\norm{\mathcal I(W)}_{H^2}\le C\norm{(\check\rho_0,\check{\bm u}_0,\check\theta_0,\check{\bm b}_0)}_{H^2}$.
\end{lemma}
\begin{proof}
For $W\in\mathcal C_R$, one has
\begin{align}\label{4.20}
 &\norm{e^{5s\Phi/6}W}_{L^\infty H^2}\le R,\qquad
   \norm{e^{5s\Phi/6}(\bm u,\theta,\bm b)}_{L^2H^3}\le R, \qquad \norm{e^{5s\Phi/6}(\partial_tW+\ub\cdot\nabla W)}_{L^2H^1}\le R,\notag\\
 &\norm{e^{5s\Phi/12}(Y-I)}_{CH^3}
       +\norm{e^{5s\Phi/12}(Z-I)}_{CH^3}\le C_sR.
\end{align}

\emph{Step 1: estimates on pressure. }
Choose $R_0$ sufficiently small that $(\rb+\rho,\tb+\theta)$ stays in a fixed compact neighborhood of $(\rb,\tb)$ contained in $(0,\infty)^2$, on which the first three derivatives of $P$ are bounded. From the Taylor formula,
\begin{align*}
 &P(\rb+\rho,\tb+\theta)-\Pb-\Prb\rho-\Ptb\theta\\
 &\quad=\int_0^1(1-\tau)\bigl[
 P_{\rho\rho}(\rb+\tau\rho,\tb+\tau\theta)\rho^2
 +2P_{\rho\theta}(\rb+\tau\rho,\tb+\tau\theta)\rho\theta
+P_{\theta\theta}(\rb+\tau\rho,\tb+\tau\theta)\theta^2\bigr]\,d\tau,
\end{align*}
and
\begin{align*}
 &\partial_j\bigl[P(\rb+\rho,\tb+\theta)-\Pb-\Prb\rho-\Ptb\theta\bigr] =\bigl[P_\rho(\rb+\rho,\tb+\theta)-\Prb\bigr]\partial_j\rho
       +\bigl[P_\theta(\rb+\rho,\tb+\theta)-\Ptb\bigr]\partial_j\theta.
\end{align*}
Differentiating once more, we also have
\begin{align*}
 &\partial_{jk}\bigl[P(\rb+\rho,\tb+\theta)-\Pb-\Prb\rho-\Ptb\theta\bigr] =\bigl[P_\rho(\rb+\rho,\tb+\theta)-\Prb\bigr]\partial_{jk}\rho
       +\bigl[P_\theta(\rb+\rho,\tb+\theta)-\Ptb\bigr]\partial_{jk}\theta\\
 &\qquad+P_{\rho\rho}(\rb+\rho,\tb+\theta)\partial_j\rho\partial_k\rho
       +P_{\theta\theta}(\rb+\rho,\tb+\theta)\partial_j\theta\partial_k\theta +P_{\rho\theta}(\rb+\rho,\tb+\theta)
                (\partial_j\rho\partial_k\theta+\partial_j\theta\partial_k\rho),
\end{align*}
which implies
\begin{align}\label{4.22}
 \norm{P(\rb+\rho,\tb+\theta)-\Pb-\Prb\rho-\Ptb\theta}_{H^2}
       &\le C\norm{(\rho,\theta)}_{H^2}^2,\notag\\
 \norm{P(\rb+\rho,\tb+\theta)-\Pb}_{H^2}
       &\le C\norm{(\rho,\theta)}_{H^2}.
\end{align}
 
\emph{Step 2: estimates on $f_\rho$.}
Using the $H^2$ algebra estimate, one has
\[
 \left\|\sum_{i,j}(\partial_iZ_j)(t,Y(t,\cdot))\partial_ju_i\right\|_{H^2}
 \le C\norm{\bm u}_{H^3},
\]
which implies
\begin{align*}
 e^{s\Phi}\left\|\rho\sum_{i,j}(\partial_iZ_j)(t,Y(t,\cdot))\partial_ju_i\right\|_{H^2}
 &\le C\norm{e^{5s\Phi/6}\rho}_{H^2}\norm{e^{5s\Phi/6}\bm u}_{H^3},\\
 e^{s\Phi}\left\|\rb\sum_{i,j}\bigl[(\partial_iZ_j)(t,Y(t,\cdot))-\delta_{ij}\bigr]\partial_ju_i\right\|_{H^2}
 &\le C\norm{e^{5s\Phi/12}(Y-I)}_{H^3}\norm{e^{5s\Phi/6}\bm u}_{H^3}.
\end{align*}
Using $L^\infty$ in time for the first factor and $L^2$ for the second, we conclude
\begin{equation}\label{4.23}
 \norm{e^{s\Phi}f_\rho}_{L^2H^2}\le C_sR^2.
\end{equation}
 
\emph{Step 3: estimates on $\bm f_{\bm u}$.}
For the time product,
\[
 e^{7s\Phi/6}\norm{\rho(\partial_t\bm u+\ub\cdot\nabla\bm u)}_{H^1}
 \le C\norm{e^{5s\Phi/6}\rho}_{H^2}
          \norm{e^{5s\Phi/6}(\partial_t\bm u+\ub\cdot\nabla\bm u)}_{H^1}.
\]
For $i,j,k =1,\ldots, d$, one has
\begin{align*}
 &e^{7s\Phi/6}\norm{[(\nabla Z_j\cdot\nabla Z_k)(t,Y(t,\cdot))-\delta_{jk}]\partial_{jk}u_i}_{H^1} \le C\norm{e^{5s\Phi/12}(Y-I)}_{H^3}\norm{e^{5s\Phi/6}\bm u}_{H^3},\\
 &e^{7s\Phi/6}\norm{(\Delta Z_k)(t,Y(t,\cdot))\partial_ku_i}_{H^1} \le C\norm{e^{5s\Phi/12}(Z-I)}_{H^3}\norm{e^{5s\Phi/6}\bm u}_{H^3}.
\end{align*}
Similarly, the two gradient-divergence terms satisfy
\begin{align*}
 &e^{7s\Phi/6}\norm{[(\partial_iZ_j)(t,Y(t,\cdot))(\partial_lZ_k)(t,Y(t,\cdot))
             -\delta_{ij}\delta_{lk}]\partial_{jk}u_l}_{H^1} \le C\norm{e^{5s\Phi/12}(Y-I)}_{H^3}\norm{e^{5s\Phi/6}\bm u}_{H^3},\\
 &e^{7s\Phi/6}\norm{(\partial_{il}Z_k)(t,Y(t,\cdot))\partial_ku_l}_{H^1} \le C\norm{e^{5s\Phi/12}(Z-I)}_{H^3}\norm{e^{5s\Phi/6}\bm u}_{H^3}.
\end{align*}
Summing finite indices and multiplying by $\mu$ or $|\lambda'+\mu|$, we arrive at the corresponding bounds in \eqref{4.11}.

It follows from \eqref{4.22} that
\begin{align*}
 &e^{7s\Phi/6}\left\|\sum_j(\partial_iZ_j)(t,Y(t,\cdot))
         \partial_j[P(\rb+\rho,\tb+\theta)-\Pb-\Prb\rho-\Ptb\theta]\right\|_{H^1} \le C\norm{e^{5s\Phi/6}(\rho,\theta)}_{H^2}^2,\\
 &e^{7s\Phi/6}\left\|\Prb\sum_j[(\partial_iZ_j)(t,Y(t,\cdot))-\delta_{ij}]\partial_j\rho\right\|_{H^1}\le C\norm{e^{5s\Phi/12}(Y-I)}_{H^3}\norm{e^{5s\Phi/6}\rho}_{H^2},\\
 &e^{7s\Phi/6}\left\|\Ptb\sum_j[(\partial_iZ_j)(t,Y(t,\cdot))-\delta_{ij}]\partial_j\theta\right\|_{H^1}\le C\norm{e^{5s\Phi/12}(Y-I)}_{H^3}\norm{e^{5s\Phi/6}\theta}_{H^2}.
\end{align*}
For the linear magnetic difference, using \eqref{4.17} with $j=1$ on the function $(\nabla\times\bm b)\times\bb$, we obtain
\begin{align*}
 e^{7s\Phi/6}\norm{(\nabla\times\bm b)(t,Y(t,\cdot))\times\bb-(\nabla\times\bm b)(t,\cdot)\times\bb}_{H^1} &\le Ce^{7s\Phi/6}\norm{Y-I}_{H^3}\norm{(\nabla\times\bm b)\times\bb}_{H^2}\\
 & \le C|\bb|\norm{e^{5s\Phi/12}(Y-I)}_{H^3}\norm{e^{5s\Phi/6}\bm b}_{H^3}.
\end{align*}
And for the last term, we have
\begin{align*}
 &e^{7s\Phi/6}\norm{(\nabla\times\bm b)(t,Y(t,\cdot))\times\bm b(t,Y(t,\cdot))}_{H^1} \le Ce^{7s\Phi/6}\norm{\nabla\times\bm b}_{H^1}\norm{\bm b}_{H^2}
 \le C\norm{e^{5s\Phi/6}\bm b}_{H^2}^2.
\end{align*}
Every time product above consists either of an $L^\infty$ factor and an $L^2$ factor, each bounded by $C_sR$, or of two $L^\infty$ factors, whose $L^2(0,T)$ product is bounded by $C_s\sqrt T R^2$. We have estimated the time term, both parts of each second-order difference, all three pressure terms and both magnetic terms. Therefore
\begin{equation}\label{4.24}
 \norm{e^{7s\Phi/6}\bm f_{\bm u}}_{L^2H^1}\le C_sR^2.
\end{equation}

\emph{Step 4:  estimates on $f_\theta$.}
By the continuous multiplication $H^2(\mathbb{T}_L)\times H^1(\mathbb{T}_L)\to H^1(\mathbb{T}_L)$, we have
\[
 e^{7s\Phi/6}\norm{c_v\rho(\partial_t\theta+\ub\cdot\nabla\theta)}_{H^1}
 \le C\norm{e^{5s\Phi/6}\rho}_{H^2}
        \norm{e^{5s\Phi/6}(\partial_t\theta+\ub\cdot\nabla\theta)}_{H^1}.
\]
The two diffusion terms satisfy, respectively,
\begin{align*}
 &e^{7s\Phi/6}\norm{[(\nabla Z_j\cdot\nabla Z_k)(t,Y(t,\cdot))-\delta_{jk}]\partial_{jk}\theta}_{H^1} \le C\norm{e^{5s\Phi/12}(Y-I)}_{H^3}\norm{e^{5s\Phi/6}\theta}_{H^3},\\
 &e^{7s\Phi/6}\norm{(\Delta Z_k)(t,Y(t,\cdot))\partial_k\theta}_{H^1} \le C\norm{e^{5s\Phi/12}(Z-I)}_{H^3}\norm{e^{5s\Phi/6}\theta}_{H^3}.
\end{align*}
And for $i,j=1,\ldots, d$, we also deduce
\begin{align*}
 &e^{7s\Phi/6}\left\|\Pb\sum_{i,j}[(\partial_iZ_j)(t,Y(t,\cdot))-\delta_{ij}]\partial_ju_i\right\|_{H^1} \le C\norm{e^{5s\Phi/12}(Y-I)}_{H^3}\norm{e^{5s\Phi/6}\bm u}_{H^2},\\
 &e^{7s\Phi/6}\left\|[P(\rb+\rho,\tb+\theta)-\Pb]
                     \sum_{i,j}(\partial_iZ_j)(t,Y(t,\cdot))\partial_ju_i\right\|_{H^1} \le C\norm{e^{5s\Phi/6}(\rho,\theta)}_{H^2}\norm{e^{5s\Phi/6}\bm u}_{H^2}.
\end{align*}

For each $i,j$, and for $r=1,2$, it holds that
\begin{align*}
 &\left\|\sum_k\bigl[(\partial_iZ_k)(t,Y(t,\cdot))\partial_ku_j
                  +(\partial_jZ_k)(t,Y(t,\cdot))\partial_ku_i\bigr]\right\|_{H^r}
       \le C\norm{\bm u}_{H^{r+1}},\\
 &\left\|\sum_{i,j}(\partial_iZ_j)(t,Y(t,\cdot))\partial_ju_i\right\|_{H^r}
       \le C\norm{\bm u}_{H^{r+1}},
\end{align*}
which implies
\begin{align*}
 &e^{7s\Phi/6}\left\|\frac\mu2\sum_{i,j}
   \left[\sum_k\bigl((\partial_iZ_k)(t,Y(t,\cdot))\partial_ku_j
                   +(\partial_jZ_k)(t,Y(t,\cdot))\partial_ku_i\bigr)\right]^2\right\|_{H^1}\le C\norm{e^{5s\Phi/6}\bm u}_{H^2}\norm{e^{5s\Phi/6}\bm u}_{H^3},\\
 &e^{7s\Phi/6}\left\|\lambda'\left[\sum_{i,j}(\partial_iZ_j)(t,Y(t,\cdot))\partial_ju_i\right]^2\right\|_{H^1}\le C\norm{e^{5s\Phi/6}\bm u}_{H^2}\norm{e^{5s\Phi/6}\bm u}_{H^3}.
\end{align*}
By the boundedness of the composition operator, Joule heating has the estimate
\[
 e^{7s\Phi/6}\norm{\nu|(\nabla\times\bm b)(t,Y(t,\cdot))|^2}_{H^1}
 \le C\norm{e^{5s\Phi/6}\bm b}_{H^2}\norm{e^{5s\Phi/6}\bm b}_{H^3},
\]
Therefore, we finally get
\begin{equation}\label{4.25}
 \norm{e^{7s\Phi/6}f_\theta}_{L^2H^1}\le C_sR^2.
\end{equation}

\emph{Step 5: estimates on $\bm f_{\bm b}$.}
Applying \eqref{4.17} in Lemma~\ref{lemma 4.2} with $j=2$, we obtain
\[
 \norm{\bm u(t,Z(t,\cdot))-\bm u(t,\cdot)}_{H^2}
 \le C\norm{Z-I}_{H^3}\norm{\bm u}_{H^3},
\]
which implies
\begin{align*}
 &e^{s\Phi}\norm{\nabla\times[(\bm u(t,Z(t,\cdot))-\bm u(t,\cdot))\times\bb]}_{H^1} \le C|\bb|\norm{e^{5s\Phi/12}(Z-I)}_{H^3}\norm{e^{5s\Phi/6}\bm u}_{H^3}.
\end{align*}
Then, by the boundedness of the composition operator, together with the algebra property of $H^2$, we have
\begin{align*}
 e^{s\Phi}\norm{\nabla\times[\bm u(t,Z(t,\cdot))\times\bm b(t,\cdot)]}_{H^1}
 &\le Ce^{s\Phi}\norm{\bm u(t,Z(t,\cdot))}_{H^2}\norm{\bm b}_{H^2} \le C\norm{e^{5s\Phi/6}\bm u}_{H^2}\norm{e^{5s\Phi/6}\bm b}_{H^2}.
\end{align*}
Consequently, we obtain
\begin{equation}\label{4.26}
 \norm{e^{s\Phi}\bm f_{\bm b}}_{L^2H^1}\le C_sR^2.
\end{equation}

Let $F$ denote the expression inside the curl in \eqref{4.13}. Then $F\in H^2(\T)^m$ and $\bm f_{\bm b}=\nabla\times F\in H^1(\T)^d$. By periodicity on \(\mathbb T_L\),
$
\int_{\mathbb T_L} \bm f_{\bm b}(t,x)\,dx=0.
$
Furthermore, $\Div \bm f_{\bm b}=0$ in the sense of distributions, since mixed distributional derivatives commute.  
Thus, we have
$
\bm f_{\bm b}(t)\in H^1_{\sigma,0}(\mathbb T_L).
$

\emph{Step 6: the initial data.}
Let \(g\) denote any one of the first three components of the initial perturbation, that is, \(g\in\{\check\rho_0,\check{\bm u}_0,\check\theta_0\}\). By Proposition~\ref{prop 4.1}, the map \(Y(0,\cdot)\) is an \(H^3\) diffeomorphism sufficiently close to the identity. Hence, according to Lemma~\ref{lemma 4.2}, the composition operator \(g\mapsto g(Y(0,\cdot))\) is bounded on \(H^2\), and therefore
\[
\norm{g(Y(0,\cdot))}_{H^2}\le C\norm{g}_{H^2}.
\]
The magnetic component remains unchanged, so its \(H^2\) norm is trivially bounded. Summing over the four components, we obtain
\[
\norm{\mathcal I(W)}_{H^2}
\le C\norm{(\check\rho_0,\check{\bm u}_0,\check\theta_0,\check{\bm b}_0)}_{H^2}.
\]
Combining this bound with estimates \eqref{4.23}--\eqref{4.26} completes the proof of the asserted inequality in the four weighted source spaces.
\end{proof}

Let \(C_G\), \(C_I\), and \(C_Q\) denote the constants associated with the linear selection, the composition-extension, and the quadratic source, respectively. It follows from \eqref{4.3} and Lemma~\ref{lemma 4.3} that
\begin{equation}\label{4.27}
\|\mathcal F(W)\|_{\mathcal B_s}\le C_G(C_I\delta+C_QR^2).
\end{equation}
To ensure that the right-hand side is bounded by \(R\), we only need to choose positive \(R\) and \(\delta\) satisfying
$
R\le\min\{R_0,(2C_GC_Q)^{-1}\}, 
\delta\le\frac{R}{2C_GC_I}.
$
Furthermore, the terminal, divergence free, and zero-mean constraints are inherited from the linear theorem. Consequently, \(\mathcal F(\mathcal C_R)\subset\mathcal C_R\).

\subsection{Continuity and existence of a fixed point}
We use the topology of $\mathcal E=L^2(Q_T)^{2+2d}$. The compactness argument is of the same type as the compressible Navier--Stokes argument in \cite{EGG2016}. However, the magnetic source and coordinate choice here are different.
\begin{proposition}\label{prop 4.4}
The set $\mathcal C_R$ is a nonempty convex compact subset of $\mathcal E$, and the map $\mathcal F:\mathcal C_R\to\mathcal C_R$ is continuous in this topology.
\end{proposition}
\begin{proof}
It is clear that $0\in\mathcal C_R$, and $\mathcal C_R$ is convex. For every $W\in\mathcal C_R$, by virtue of \eqref{4.1}, we have
\[
\norm W_{L^\infty H^2}+\norm{\partial_tW}_{L^2H^1}\le CR.
\]
Hence, for $0<\alpha<1$ and $t,t'\in[0,T]$,
\[
\norm{W(t')-W(t)}_{H^{2-\alpha}}
\le C\norm{W(t')-W(t)}_{H^1}^{\alpha}\norm{W(t')-W(t)}_{H^2}^{1-\alpha}
\le CR|t'-t|^{\alpha/2}.
\]
Since $H^2\subset\subset H^{2-\alpha}$, the Arzel\`a--Ascoli theorem yields relative compactness in $CH^{2-\alpha}$, and thus in $\mathcal E$. If $W_n\in\mathcal C_R$ and $W_n\to W$ in $\mathcal E$, weak and weak-star compactness of the weighted fields, followed by identification on each interval $[0,T-\tau]$, $\tau>0$, gives
\[
\norm W_{\mathcal B_s}\le\liminf_{n\to\infty}\norm{W_n}_{\mathcal B_s}\le R.
\]
Moreover, convergence in $CH^{2-\alpha}$ and distributional testing imply
$
W(T)=0, \Div\bm b=0, \int_\T\bm b(t,x)\,dx=0.
$
Thus $\mathcal C_R$ is closed in $\mathcal E$ and therefore compact.

Let now $W_n\to W$ in $\mathcal E$, with $W_n,W\in\mathcal C_R$, and fix $0<\alpha<1/2$. Compactness and interpolation give
\begin{equation}\label{4.28}
W_n\longrightarrow W\quad\text{in }CH^{2-\alpha},\qquad
(\bm u_n,\theta_n,\bm b_n)\longrightarrow(\bm u,\theta,\bm b)\quad\text{in }L^2H^{3-\alpha}.
\end{equation}
In particular, $W_n\to W$ in $CL^\infty$. Denote by $Y_n,Z_n$ the maps associated with $W_n$. It follows from \eqref{4.4} that
\[
\norm{Y_n-Y}_{CH^{3-\alpha}}\le\sqrt T\norm{\bm u_n-\bm u}_{L^2H^{3-\alpha}}\longrightarrow0.
\]
The embedding $H^{3-\alpha}\hookrightarrow W^{1,\infty}$, Proposition~\ref{prop 4.1}, and the identities
$
Y_n(t,Z_n(t,x))=x, D_xZ_n(t,x)=[(D_xY_n)(t,Z_n(t,x))]^{-1}
$
then yield $Z_n\to Z$ in $CW^{1,\infty}$. Set 
$ 
A_{n,ij}(t,x)=(\partial_iZ_{n,j})(t,Y_n(t,x)),  A_{ij}(t,x)=(\partial_iZ_j)(t,Y(t,x)).
$
By \eqref{4.5} and continuity of matrix inversion in $H^{2-\alpha}$,
\[
A_{n,ij}\longrightarrow A_{ij}\quad\text{in }CH^{2-\alpha},\qquad
\partial_kA_{n,ij}\longrightarrow\partial_kA_{ij}\quad\text{in }CL^2.
\]
Furthermore,
\begin{align*}
B_{n,il,k}(t,x)&:=(\partial_{il}Z_{n,k})(t,Y_n(t,x))=\sum_jA_{n,ij}\partial_jA_{n,lk}\\
&\longrightarrow\sum_jA_{ij}\partial_jA_{lk}=(\partial_{il}Z_k)(t,Y(t,x))=:B_{il,k}(t,x)\quad\text{in }CL^2.
\end{align*}
For $g\in H^2$, choose $g^{(k)}\in C^\infty(\T)$ such that $g^{(k)}\to g$ in $H^2$. Lemma~\ref{lemma 4.2} and the chain rule imply
\begin{align*}
\norm{g(Y_n(0,\cdot))-g(Y(0,\cdot))}_{H^2}
&\le C\norm{g-g^{(k)}}_{H^2}+\norm{g^{(k)}(Y_n(0,\cdot))-g^{(k)}(Y(0,\cdot))}_{H^2}\\
&\le C\norm{g-g^{(k)}}_{H^2}+C_k\norm{Y_n(0,\cdot)-Y(0,\cdot)}_{H^2}.
\end{align*}
Letting first $n\to\infty$ and then $k\to\infty$, and applying the result to the components in \eqref{4.14}, we obtain
\begin{equation}\label{4.29}
\mathcal I(W_n)\longrightarrow\mathcal I(W)\quad\text{in }H^2.
\end{equation}
We next identify the limit of the sources in \eqref{4.10}--\eqref{4.13}. For a scalar component $v_n=u_{n,i}$ or $v_n=\theta_n$, the corresponding limit $v$ satisfies
\[
(\partial_t+\ub\cdot\nabla)v_n\rightharpoonup(\partial_t+\ub\cdot\nabla)v\quad\text{in }L^2H^1.
\]
Thus, for every $\psi\in C_c^\infty(Q_T)$,
\begin{align*}
&\int_{Q_T}\big[\rho_n(\partial_t+\ub\cdot\nabla)v_n-\rho(\partial_t+\ub\cdot\nabla)v\big]\psi\\
&\quad=\int_{Q_T}(\rho_n-\rho)(\partial_t+\ub\cdot\nabla)v_n\psi+\int_{Q_T}(\partial_t+\ub\cdot\nabla)(v_n-v)\rho\psi\longrightarrow0.
\end{align*}
On the other hand, it follows from \eqref{4.28} and the coefficient limits above, one has
\begin{align*}
A_{n,ij}A_{n,lk}\partial_{jk}v_n&\longrightarrow A_{ij}A_{lk}\partial_{jk}v\quad\text{in }L^2L^2, ~~
B_{n,il,k}\partial_kv_n\longrightarrow B_{il,k}\partial_kv\quad\text{in }L^1(Q_T),\\
P(\rb+\rho_n,\tb+\theta_n)&\longrightarrow P(\rb+\rho,\tb+\theta)\quad\text{in }CL^\infty,~~
\nabla P(\rb+\rho_n,\tb+\theta_n)\longrightarrow\nabla P(\rb+\rho,\tb+\theta)\quad\text{in }CL^2.
\end{align*}
These limits identify the diffusion and pressure terms, as well as $f_\rho(W_n)$, in distributions. For completeness, if $g_n\to g$ in $L^2(Q_T)$ and $V_n\to V$ uniformly, where the maps $V_n(t,\cdot)$ and $V(t,\cdot)$ have uniformly bounded inverse Jacobians, then, for every smooth $g^{(k)}$,
\begin{align*}
&\norm{g_n(t,V_n(t,x))-g(t,V(t,x))}_{L^2(Q_T)}\\
&\quad\le C\norm{g_n-g}_{L^2(Q_T)}+C\norm{g-g^{(k)}}_{L^2(Q_T)}
+\norm{g^{(k)}(t,V_n(t,x))-g^{(k)}(t,V(t,x))}_{L^2(Q_T)}.
\end{align*}
Smooth approximation, with $n\to\infty$ preceding $k\to\infty$, therefore gives convergence to zero. Applying this observation with $V_n=Y_n$ or $V_n=Z_n$, we find
\begin{align*}
(\nabla\times\bm b_n)(t,Y_n(t,x))&\longrightarrow(\nabla\times\bm b)(t,Y(t,x))\quad\text{in }L^2(Q_T),\\
\bm b_n(t,Y_n(t,x))&\longrightarrow\bm b(t,Y(t,x)),\quad
\bm u_n(t,Z_n(t,x))\longrightarrow\bm u(t,Z(t,x))\quad\text{in }L^2(Q_T).
\end{align*}
In addition, $A_{n,ik}\partial_k u_{n,j}\to A_{ik}\partial_k u_j$ in $L^2(Q_T)$. Consequently, the Lorentz terms converge in $L^1(Q_T)$, and so do all the quadratic heating terms, by
\[
\norm{|h_n|^2-|h|^2}_{L^1(Q_T)}\le\big(\norm{h_n}_{L^2(Q_T)}+\norm h_{L^2(Q_T)}\big)\norm{h_n-h}_{L^2(Q_T)}.
\]
Finally, letting $\bm b_n\to\bm b$ in $CL^\infty$ yields
\begin{align*}
&[\bm u_n(t,Z_n(t,x))-\bm u_n(t,x)]\times\bb+\bm u_n(t,Z_n(t,x))\times\bm b_n(t,x)\\
&\qquad\longrightarrow[\bm u(t,Z(t,x))-\bm u(t,x)]\times\bb+\bm u(t,Z(t,x))\times\bm b(t,x)\quad\text{in }L^2(Q_T).
\end{align*}
Taking the curl in distributions identifies the limit in \eqref{4.13}. Together, these relations show
\[
f(W_n)\longrightarrow f(W)\quad\text{in }\mathcal D'(Q_T).
\]
By Lemma~\ref{lemma 4.3}, the sequence $f(W_n)$ is bounded in the reflexive space $\mathcal Y_s$. 
Every subsequence has a further weakly convergent subsequence, whose limit is the distributional limit just identified. Hence the whole sequence satisfies
\[
f(W_n)\rightharpoonup f(W)\quad\text{in }\mathcal Y_s.
\]
It follows from the bounded linearity of $\mathcal G_{\rm st}$, \eqref{4.15}, and \eqref{4.29} that
\[
\mathcal F(W_n)=\mathcal G_{\rm st}(\mathcal I(W_n),f(W_n))\rightharpoonup\mathcal G_{\rm st}(\mathcal I(W),f(W))=\mathcal F(W)\quad\text{in }\mathcal E.
\]
Since $\mathcal F(W_n)\in\mathcal C_R$ and $\mathcal C_R$ is compact in $\mathcal E$, every subsequence has a further strongly convergent subsequence, whose limit must be $\mathcal F(W)$ by the preceding weak convergence. It follows that $\norm{\mathcal F(W_n)-\mathcal F(W)}_{\mathcal E}\to0$.
\end{proof}

\begin{proof}[Proof of Theorem~\ref{thm 1.1}]
Choose \(R\) and \(\delta\) as in \eqref{4.27}. By Proposition~\ref{prop 4.4} and Schauder's fixed point theorem, there exists a fixed point \(W\in\mathcal C_R\). Applying Theorem~\ref{thm 3.7}, we select the associated controls such that \(W\) solves \eqref{4.9}  with sources \eqref{4.10}--\eqref{4.13} and initial data \eqref{4.14}.

Now recover the physical fields via \eqref{4.7}. The first three controls in the extended primitive equations read
\[
 \check v_\rho(t,x)=\chi(Z(t,x))v_\rho(t,Z(t,x)),~
 \check{\bm v}_{\bm u}(t,x)=\chi(Z(t,x))\bm v_{\bm u}(t,Z(t,x)),~
 \check v_\theta(t,x)=\chi(Z(t,x))v_\theta(t,Z(t,x)),
\]
whereas the magnetic control remains \(\nabla\times(\chi v_{\bm b})\) in physical coordinates. Due to Proposition~\ref{prop 4.1}, these sources vanish in a neighborhood of \(\overline\Omega\). 
Substituting the explicit definitions of the source terms \(f_\rho\), \(\bm f_{\bm u}\), and \(f_\theta\) from \eqref{4.10}--\eqref{4.12} into the first three equations of \eqref{4.9}, we recover the first three equations of \eqref{2.2}.
As for the fourth, 
\[
 \nabla\times(\bm u\times\bb)+\bm f_{\bm b}(W)
 =\nabla\times[\bm u(t,Z(t,x))\times\bb]
   +\nabla\times[\bm u(t,Z(t,x))\times\bm b(t,x)].
\]
This is precisely the induction equation for \(\check{\bm u}(t,x)=\bm u(t,Z(t,x))\) and \(\check{\bm b}(t,x)=\bm b(t,x)\), whose constraint is inherited without any spatial change of the magnetic field.

The identity \(Y(0,Z(0,x))=x\) together with \eqref{4.14} yields the prescribed initial data. Since \(Y(T,x)=Z(T,x)=x\) and \(W(T)=0\), the physical terminal state coincides with the equilibrium. At this fixed point, Theorem~\ref{thm 3.7} and Remark~\ref{remark 3.1} enhance the time regularity to \(CH^2\), while the bounds \(Y,Z\in CH^3\) preserve the spatial norms. In addition,
\[
 \partial_tY(t,x)=\bm u(t,x)-\ub\cdot\nabla(Y(t,x)-x),\qquad
 \partial_tZ(t,x)=-D_xZ(t,x)(\partial_tY)(t,Z(t,x))
\]
are bounded in \(L^\infty H^2\). 
By \eqref{4.16}, it follows that the reconstructed fields preserve \(H^1H^1\) regularity, and hence \eqref{1.5} is obtained.

Finally, restriction to \(\Omega\) and the trace theorem provide the controls and regularities stated in Theorem~\ref{thm 1.1}, with the density trace used only on \(\Gamma_-(\bm U)\). The magnetic field itself is the solenoidal lifting required in Remark~\ref{remark 1.3}, and its separate boundary fluxes vanish by Remark~\ref{remark 1.2}. The endpoint compatibility follows from \eqref{1.3}--\eqref{1.5}. The small $L^\infty H^2$ bound, Sobolev embedding, and the choice of $R_0$ give $\widetilde\rho,\widetilde\theta>0$.
\end{proof}

\section*{Acknowledgement}
The author is partially supported by the National Natural Science Foundation of China (Grant No. 12471423), by the Guangdong Basic and Applied Basic
Research Foundation (Grant No. 2025A1515010485), by the Shenzhen Science and Technology (Program Nos. JCYJ20241202124209011 and
RCJC20221008092753082), and by Shenzhen University 2035 Program for Excellent Research (Grant No. 2024C008).


\vskip8mm


\begin{thebibliography}{99}

\bibitem{A2011} E. V. Amosova, \emph{Exact local controllability for equations of the dynamics of a viscous gas}, (Russian) Differ. Uravn., 47 (2011), pp. 1754--1772.

\bibitem{B2014} M. Badra, \emph{Local controllability to trajectories of the magnetohydrodynamic equations}, J. Math. Fluid Mech., 16 (2014), pp. 631--660.

\bibitem{BEG2016} M. Badra, S. Ervedoza and S. Guerrero, \emph{Local controllability to trajectories for non-
homogeneous incompressible Navier-Stokes equations}, Ann. Inst. H. Poincar\'{e} Non Linear Anal.,
33 (2016), pp. 529--574.

\bibitem{BHPS2003} V. Barbu, T. Hav\^{a}rneanu, C. Popa and S. S. Sritharan, \emph{Exact controllability for the magnetohydrodynamic equations}, Comm. Pure Appl. Math., 56 (2003), pp. 732--783.

\bibitem{BHPS2005} V. Barbu, T. Hav\^{a}rneanu, C. Popa and S. S. Sritharan, \emph{Local exact controllability for the magnetohydrodynamic equations revisited}, Adv. Differential Equations, 10 (2005), pp. 481--504.

\bibitem{C1970} H. Cabannes, \emph{Theoretical Magnetofluiddynamics}, Academic Press, New York, London, 1970.

\bibitem{CM2015} S. Chowdhury and D. Mitra, \emph{Null controllability of the linearized compressible Navier--Stokes equations using moment method},
J. Evol. Equ., 15 (2015), pp. 331--360.

\bibitem{CMRR2014} S. Chowdhury, D. Mitra, M. Ramaswamy and M. Renardy, \emph{Null controllability of the linearized compressible Navier-Stokes system in one dimension}, J. Differential Equations, 257 (2014), pp. 3813--3849.

\bibitem{CRR2012} S. Chowdhury, M. Ramaswamy and J. P. Raymond, \emph{Controllability and stabilizability of the linearized compressible Navier-Stokes system in one dimension}, SIAM J. Control Optim., 50 (2012), pp. 2959--2987.

\bibitem{CM2010} M. Costabel and A. McIntosh,
\emph{On Bogovski\u\i\ and regularized Poincar\'e integral operators for de Rham complexes on Lipschitz domains}, Math. Z., 265 (2010), pp. 297--320.

\bibitem{EGG2016} S. Ervedoza, O. Glass and S. Guerrero, \emph{Local exact controllability for the two- and three-dimensional compressible Navier-Stokes equations}, Comm. Partial Differential Equations, 41 (2016), pp. 1660--1691.

\bibitem{EGGP2012} S. Ervedoza, O. Glass, S. Guerrero and J. P. Puel, \emph{Local exact controllability for the one-dimensional compressible Navier-Stokes equation}, Arch. Rational Mech. Anal., 206 (2012), pp. 189--238.

\bibitem{ES2018} S. Ervedoza and M. Savel, \emph{Local boundary controllability to trajectories for the 1D compressible Navier--Stokes equations}, ESAIM Control Optim. Calc. Var., 24 (2018), pp. 211--235.

\bibitem{G2013} I. Galan, \emph{Approximate controllability of the magnetohydrodynamic equations on the three-dimensional torus}, J. Optim. Theory Appl., 159 (2013), pp. 231--245.

\bibitem{G2007} O. Glass, \emph{On the controllability of the 1-D isentropic Euler equation}, J. Eur. Math. Soc.,
9 (2007), pp. 427--486.

\bibitem{G2014} O. Glass, \emph{On the controllability of the non-isentropic 1-D Euler equation}, J. Differential Equations,
257(2014), pp. 638--719.

\bibitem{HPS2006} T. Hav\^{a}rneanu, C. Popa and S. S. Sritharan, \emph{Exact internal controllability for the magnetohydrodynamic equations in multi-connected domains}, Adv. Differential Equations, 11 (2006), pp. 893--929.

\bibitem{HPS2007} T. Hav\^{a}rneanu, C. Popa and S. S. Sritharan, \emph{Exact internal controllability for the two-dimensional magnetohydrodynamic equations}, SIAM J. Control Optim., 46 (2007), pp. 1802--1830.

\bibitem{KNV2022} I. Kukavica, M. Novack and V. Vicol, \emph{Exact boundary controllability for the ideal magneto-hydrodynamic equations}, J. Differential Equations, 318 (2022), pp. 94--112.

\bibitem{KL1965} A. G. Kulikovskiy and G. A. Lyubimov, \emph{Magnetohydrodynamics}, Addison-Wesley, Reading,
Massachusetts, 1965.

\bibitem{KSV2025} K. Koike, F. Sueur and G. Vergara-Hermosilla, \emph{Local exact Lagrangian controllability for the one-dimensional barotropic compressible Navier--Stokes equations}, SIAM J. Control Optim., 63 (2025), pp. 2177--2191.

\bibitem{Lions1998} P. L. Lions, \emph{Mathematical Topics in Fluid Mechanics, Vol.~2: Compressible Models}.
Oxford University Press, New York, 1998.

\bibitem{LR2003} T. T. Li and B. P. Rao, \emph{Exact boundary controllability for quasi-linear hyperbolic systems},
SIAM J. Control Optim., 41 (2003), pp. 1748--1755.

\bibitem{Maity2015} D. Maity, \emph{Some controllability results for linearized compressible Navier-Stokes system},
ESAIM Control Optim. Calc. Var., 21 (2015), pp. 1002--1028.

\bibitem{M2019} N. Molina, \emph{Local exact boundary controllability for the compressible Navier-Stokes equations}, SIAM J. Control Optim., 57 (2019), pp. 2152--2184.

\bibitem{N2011} H. Nersisyan, \emph{Controllability of the 3D compressible Euler system}, Comm. Partial Differential Equations, 36 (2011), pp. 1544--1564.

\bibitem{R2025} M. Rissel, \emph{Exact controllability of incompressible ideal magnetohydrodynamics in two dimensions}, Pure Appl. Anal., 7 (2025), pp. 661--700.

\bibitem{RW2021} M. Rissel and Y.-G. Wang, \emph{Global exact controllability of ideal incompressible magnetohydrodynamic flows through a planar duct}, ESAIM Control Optim. Calc. Var., 27 (2021), Paper No. 103, 24 pp.

\bibitem{RW2024} M. Rissel and Y.-G. Wang, \emph{Small-time global approximate controllability for incompressible MHD with coupled Navier slip boundary conditions}, J. Math. Pures Appl., 190 (2024), Paper No. 103601.

\bibitem{T2018} Q. Tao, \emph{Local exact controllability for the planar compressible magnetohydrodynamic equations}, SIAM J. Control Optim., 56 (2018), pp. 4461--4487.

\bibitem{TYY2025} Q. Tao, Z.-A. Yao and X. Yin, \emph{Local exact controllability to the trajectories for the two-dimensional magnetohydrodynamic system with controls acting only on the velocity field}, J. Differential Equations, 433 (2025), Paper No. 113237.

\end{thebibliography}
\end{document}